\documentclass{lmcs}
\pdfoutput=1

\usepackage{lastpage}
\lmcsdoi{22}{2}{35}
\lmcsheading{}{\pageref{LastPage}}{}{}%
{Sep.~23,~2024}{Jun.~29,~2026}{}

\usepackage[numbers]{natbib}
\usepackage{longtable,booktabs,array}
\usepackage{amssymb,amsmath}
\usepackage{stmaryrd}
\usepackage[utf8]{inputenc}
\usepackage[all]{xy}

\usepackage{tikz-cd}
\usepackage{graphicx}
\usepackage{color}

\newenvironment{theorem}{\begin{thm}}{\end{thm}}
\newenvironment{corollary}{\begin{cor}}{\end{cor}}
\newenvironment{lemma}{\begin{lem}}{\end{lem}}
\newenvironment{proposition}{\begin{prop}}{\end{prop}}

\newenvironment{definition}{\begin{defi}}{\end{defi}}

\usepackage{hyperref}

\usepackage{agdahen}

\providecommand{\tightlist}{%
  \setlength{\itemsep}{0pt}\setlength{\parskip}{0pt}}

\usepackage{esint}
\DeclareUnicodeCharacter{2192}{\ensuremath{\rightarrow}} 
\DeclareUnicodeCharacter{2190}{\ensuremath{\leftarrow}} 
\DeclareUnicodeCharacter{22EF}{\ensuremath{\cdots}} 
\DeclareUnicodeCharacter{2208}{\ensuremath{\in}} 
\DeclareUnicodeCharacter{220B}{\ensuremath{\ni}} 
\DeclareUnicodeCharacter{2200}{\ensuremath{\forall}} 
\DeclareUnicodeCharacter{2203}{\ensuremath{\exists}} 
\DeclareUnicodeCharacter{03B1}{\ensuremath{\alpha}} 
\DeclareUnicodeCharacter{03B2}{\ensuremath{\beta}} 
\DeclareUnicodeCharacter{03BE}{\ensuremath{\xi}} 
\DeclareUnicodeCharacter{03B4}{\ensuremath{\delta}} 
\DeclareUnicodeCharacter{03B5}{\ensuremath{\epsilon}} 
\DeclareUnicodeCharacter{03C6}{\ensuremath{\phi}} 
\DeclareUnicodeCharacter{03B3}{\ensuremath{\gamma}} 
\DeclareUnicodeCharacter{03B8}{\ensuremath{\theta}} 
\DeclareUnicodeCharacter{03B9}{\ensuremath{\iota}} 
\DeclareUnicodeCharacter{03BA}{\ensuremath{\kappa}} 
\DeclareUnicodeCharacter{03BB}{\ensuremath{\lambda}} 
\DeclareUnicodeCharacter{03BC}{\ensuremath{\mu}} 
\DeclareUnicodeCharacter{03BD}{\ensuremath{\nu}} 
\DeclareUnicodeCharacter{03BF}{\ensuremath{\omicron}} 
\DeclareUnicodeCharacter{03C0}{\ensuremath{\pi}} 
\DeclareUnicodeCharacter{03C8}{\ensuremath{\psi}} 
\DeclareUnicodeCharacter{03C1}{\ensuremath{\rho}} 
\DeclareUnicodeCharacter{03C3}{\ensuremath{\sigma}} 
\DeclareUnicodeCharacter{03C4}{\ensuremath{\tau}} 
\DeclareUnicodeCharacter{03C5}{\ensuremath{\upsilon}} 
\DeclareUnicodeCharacter{03C9}{\ensuremath{\omega}} 
\DeclareUnicodeCharacter{03C7}{\ensuremath{\chi}} 
\DeclareUnicodeCharacter{03B7}{\ensuremath{\eta}} 
\DeclareUnicodeCharacter{03B6}{\ensuremath{\zeta}} 
\DeclareUnicodeCharacter{0391}{\ensuremath{\Alpha}} 
\DeclareUnicodeCharacter{0392}{\ensuremath{\Beta}} 
\DeclareUnicodeCharacter{039E}{\ensuremath{\Xi}} 
\DeclareUnicodeCharacter{0394}{\ensuremath{\Delta}} 
\DeclareUnicodeCharacter{0395}{\ensuremath{\Epsilon}} 
\DeclareUnicodeCharacter{03A6}{\ensuremath{\Phi}} 
\DeclareUnicodeCharacter{0393}{\ensuremath{\Gamma}} 
\DeclareUnicodeCharacter{0398}{\ensuremath{\Theta}} 
\DeclareUnicodeCharacter{0399}{\ensuremath{\Iota}} 
\DeclareUnicodeCharacter{039A}{\ensuremath{\Kappa}} 
\DeclareUnicodeCharacter{039B}{\ensuremath{\Lambda}} 
\DeclareUnicodeCharacter{039C}{\ensuremath{\Mu}} 
\DeclareUnicodeCharacter{039D}{\ensuremath{\Nu}} 
\DeclareUnicodeCharacter{039F}{\ensuremath{\Omicron}} 
\DeclareUnicodeCharacter{03A0}{\ensuremath{\Pi}} 
\DeclareUnicodeCharacter{03A8}{\ensuremath{\Psi}} 
\DeclareUnicodeCharacter{03A1}{\ensuremath{\Rho}} 
\DeclareUnicodeCharacter{03A3}{\ensuremath{\Sigma}} 
\DeclareUnicodeCharacter{03A4}{\ensuremath{\Tau}} 
\DeclareUnicodeCharacter{03A5}{\ensuremath{\Upsilon}} 
\DeclareUnicodeCharacter{03A9}{\ensuremath{\Omega}} 
\DeclareUnicodeCharacter{03A7}{\ensuremath{\Chi}} 
\DeclareUnicodeCharacter{0397}{\ensuremath{\Eta}} 
\DeclareUnicodeCharacter{0396}{\ensuremath{\Zeta}} 
\DeclareUnicodeCharacter{2219}{\ensuremath{\bullet}} 
\DeclareUnicodeCharacter{263A}{:-)} 
\DeclareUnicodeCharacter{2193}{\ensuremath{\downarrow}} 
\DeclareUnicodeCharacter{1D48}{\ensuremath{^d}} 
\DeclareUnicodeCharacter{2153}{\ensuremath{\frac{1}{3}}} 
\DeclareUnicodeCharacter{2211}{\ensuremath{\sum}} 
\DeclareUnicodeCharacter{220F}{\ensuremath{\prod}} 
\DeclareUnicodeCharacter{2243}{\ensuremath{\simeq}} 
\DeclareUnicodeCharacter{2248}{\ensuremath{\approx}} 
\DeclareUnicodeCharacter{207B}{\ensuremath{^-}} 
\DeclareUnicodeCharacter{00B9}{\ensuremath{^1}} 
\DeclareUnicodeCharacter{2080}{\ensuremath{_0}} 
\DeclareUnicodeCharacter{2245}{\ensuremath{\cong}} 
\DeclareUnicodeCharacter{2218}{\ensuremath{\circ}} 
\DeclareUnicodeCharacter{2254}{\ensuremath{:=}} 
\DeclareUnicodeCharacter{00B7}{\ensuremath{\cdot}} 
\DeclareUnicodeCharacter{2081}{\ensuremath{_1}} 
\DeclareUnicodeCharacter{2261}{\ensuremath{\equiv}} 
\DeclareUnicodeCharacter{2227}{\ensuremath{\wedge}} 
\DeclareUnicodeCharacter{2228}{\ensuremath{\vee}} 
\DeclareUnicodeCharacter{00D7}{\ensuremath{\times}} 
\DeclareUnicodeCharacter{2205}{\ensuremath{\emptyset}} 
\DeclareUnicodeCharacter{2194}{\ensuremath{\leftrightarrow}} 
\DeclareUnicodeCharacter{2260}{\ensuremath{\neq}} 
\DeclareUnicodeCharacter{2264}{\ensuremath{\le}} 
\DeclareUnicodeCharacter{2265}{\ensuremath{\ge}} 
\DeclareUnicodeCharacter{2423}{\,} 
\DeclareUnicodeCharacter{211D}{\ensuremath{\mathbb{R}}} 
\DeclareUnicodeCharacter{266F}{\ensuremath{\sharp}} 
\DeclareUnicodeCharacter{2212}{\ensuremath{-}} 
\DeclareUnicodeCharacter{22C0}{\ensuremath{\bigwedge}} 
\DeclareUnicodeCharacter{2295}{\ensuremath{\oplus}} 
\DeclareUnicodeCharacter{2297}{\ensuremath{\otimes}} 
\DeclareUnicodeCharacter{2299}{\ensuremath{\odot}} 
\DeclareUnicodeCharacter{2102}{\ensuremath{\mathbb{C}}} 
\DeclareUnicodeCharacter{2112}{\ensuremath{\mathcal{L}}} 
\DeclareUnicodeCharacter{2130}{\ensuremath{\mathcal{E}}} 
\DeclareUnicodeCharacter{2217}{\ensuremath{*}} 
\DeclareUnicodeCharacter{2131}{\ensuremath{\mathcal{F}}} 
\DeclareUnicodeCharacter{22A8}{\ensuremath{\Vdash}} 
\DeclareUnicodeCharacter{21D2}{\ensuremath{\Rightarrow}} 
\DeclareUnicodeCharacter{2115}{\ensuremath{\mathbb{N}}} 
\DeclareUnicodeCharacter{2013}{--} 
\DeclareUnicodeCharacter{2010}{-} 
\DeclareUnicodeCharacter{2124}{\ensuremath{\mathbb{Z}}} 
\DeclareUnicodeCharacter{207F}{\ensuremath{^n}} 
\DeclareUnicodeCharacter{2099}{\ensuremath{_n}} 
\DeclareUnicodeCharacter{2032}{\ensuremath{'}} 
\DeclareUnicodeCharacter{221A}{\ensuremath{\sqrt}} 
\DeclareUnicodeCharacter{2071}{\ensuremath{^i}} 
\DeclareUnicodeCharacter{211A}{\ensuremath{\mathbb{Q}}} 
\DeclareUnicodeCharacter{21D4}{\ensuremath{\Leftrightarrow}} 
\DeclareUnicodeCharacter{00BE}{\ensuremath{\frac 34}} 
\DeclareUnicodeCharacter{2329}{\ensuremath{\langle}} 
\DeclareUnicodeCharacter{232A}{\ensuremath{\rangle}} 
\DeclareUnicodeCharacter{2286}{\ensuremath{\subseteq}} 
\DeclareUnicodeCharacter{22A2}{\ensuremath{\vdash}} 
\DeclareUnicodeCharacter{212C}{\ensuremath{\mathcal{B}}} 
\DeclareUnicodeCharacter{227A}{\ensuremath{\prec}} 
\DeclareUnicodeCharacter{22A4}{\ensuremath{\top}} 
\DeclareUnicodeCharacter{2202}{\ensuremath{\partial}} 
\DeclareUnicodeCharacter{211B}{\ensuremath{\mathcal{R}}} 
\DeclareUnicodeCharacter{222A}{\ensuremath{\cup}} 
\DeclareUnicodeCharacter{2082}{\ensuremath{_2}} 
\DeclareUnicodeCharacter{208A}{\ensuremath{_+}} 
\DeclareUnicodeCharacter{22C3}{\ensuremath{\bigcup}} 
\DeclareUnicodeCharacter{2229}{\ensuremath{\cap}} 
\DeclareUnicodeCharacter{221E}{\ensuremath{\infty}} 
\DeclareUnicodeCharacter{2074}{\ensuremath{^4}} 
\DeclareUnicodeCharacter{2084}{\ensuremath{_4}} 
\DeclareUnicodeCharacter{00B2}{\ensuremath{^2}} 
\DeclareUnicodeCharacter{00B3}{\ensuremath{^3}} 
\DeclareUnicodeCharacter{2119}{\ensuremath{\mathbb P}} 
\DeclareUnicodeCharacter{2207}{\ensuremath{\nabla}} 
\DeclareUnicodeCharacter{212F}{\ensuremath{e}} 
\DeclareUnicodeCharacter{22C1}{\ensuremath{\bigvee}} 
\DeclareUnicodeCharacter{207A}{\ensuremath{^+}} 
\DeclareUnicodeCharacter{220E}{\qed \vspace{3mm}} 
\DeclareUnicodeCharacter{22A5}{\ensuremath{\bot}} 
\DeclareUnicodeCharacter{2133}{\ensuremath{\mathcal{M}}} 
\DeclareUnicodeCharacter{2500}{\ensuremath{{-\mkern-3.4mu}}} 
\DeclareUnicodeCharacter{2110}{\ensuremath{\mathcal{I}}} 
\DeclareUnicodeCharacter{2282}{\ensuremath{\subset}} 
\DeclareUnicodeCharacter{2283}{\ensuremath{\supset}} 
\DeclareUnicodeCharacter{2070}{\ensuremath{^0}} 
\DeclareUnicodeCharacter{2502}{|} 
\DeclareUnicodeCharacter{2083}{\ensuremath{_3}} 
\DeclareUnicodeCharacter{016B}{\={u}} 
\DeclareUnicodeCharacter{22B2}{\ensuremath{\triangleleft}} 
\DeclareUnicodeCharacter{22B3}{\ensuremath{\triangleright}} 
\DeclareUnicodeCharacter{21AA}{\ensuremath{\hookrightarrow}} 
\DeclareUnicodeCharacter{2225}{\ensuremath{\|}} 
\DeclareUnicodeCharacter{2665}{\ensuremath{\heartsuit}} 
\DeclareUnicodeCharacter{208B}{\ensuremath{{_-}}} 
\DeclareUnicodeCharacter{00A6}{\ensuremath{\|}} 
\DeclareUnicodeCharacter{2571}{\ensuremath{\diagup}} 
\DeclareUnicodeCharacter{2572}{\ensuremath{\diagdown}} 
\DeclareUnicodeCharacter{2573}{X} 
\DeclareUnicodeCharacter{2514}{+} 
\DeclareUnicodeCharacter{2504}{\ensuremath{\cdots}} 
\DeclareUnicodeCharacter{2518}{+} 
\DeclareUnicodeCharacter{00AC}{\ensuremath{\neg}} 
\DeclareUnicodeCharacter{2191}{\ensuremath{\uparrow}} 
\DeclareUnicodeCharacter{223C}{\ensuremath{\sim}} 
\DeclareUnicodeCharacter{21A9}{\ensuremath{\hookleftarrow}} 
\DeclareUnicodeCharacter{2012}{--} 
\DeclareUnicodeCharacter{21A0}{\ensuremath{\twoheadrightarrow}} 
\DeclareUnicodeCharacter{222B}{\ensuremath{\int}} 
\DeclareUnicodeCharacter{21C7}{\ensuremath{\overleftarrow{\leftarrow}}} 
\DeclareUnicodeCharacter{21C9}{\ensuremath{\overrightarrow{\rightarrow}}} 
\DeclareUnicodeCharacter{21CA}{\ensuremath{\downarrow\downarrow}} 
\DeclareUnicodeCharacter{25A1}{\ensuremath{\Box}} 
\DeclareUnicodeCharacter{01B2}{\ensuremath{\mathcal{V}}} 
\DeclareUnicodeCharacter{2148}{\ensuremath{\mathbb{i}}} 
\DeclareUnicodeCharacter{231F}{\ensuremath{\lrcorner}} 
\DeclareUnicodeCharacter{226C}{\ensuremath{\between}} 
\DeclareUnicodeCharacter{2118}{\ensuremath{\mathcal P}} 
\DeclareUnicodeCharacter{0192}{\ensuremath{\fint}} 
\DeclareUnicodeCharacter{043B}{\ensuremath{\el}} 
\DeclareUnicodeCharacter{2294}{\ensuremath{\sqcup}} 
\DeclareUnicodeCharacter{2A06}{\ensuremath{\bigsqcup}} 
\DeclareUnicodeCharacter{219D}{\ensuremath{\leadsto}} 
\DeclareUnicodeCharacter{2237}{\ensuremath{{::}}} 
\DeclareUnicodeCharacter{231D}{\ensuremath{\urcorner}} 
\DeclareUnicodeCharacter{231C}{\ensuremath{\ulcorner}} 
\DeclareUnicodeCharacter{2022}{\ensuremath{\bullet}} 
\DeclareUnicodeCharacter{21A6}{\ensuremath{\mapsto}} 
\DeclareUnicodeCharacter{21A0}{\ensuremath{\twoheadrightarrow}} 
\DeclareUnicodeCharacter{25C1}{\ensuremath{\triangleleft}} 
\DeclareUnicodeCharacter{27E6}{\ensuremath{\llbracket}} 
\DeclareUnicodeCharacter{27E7}{\ensuremath{\rrbracket}} 

\newcommand{\Set}{\ensuremath{\operatorname{Set}}}
\newcommand{\Type}{\ensuremath{\operatorname{Type}}}

\newcommand{\inl}{\ensuremath{\operatorname{inl}}} 
\newcommand{\Nat}{ℕ}

\newcommand{\isntruncmap}[1]{\ensuremath{\operatorname{is\mbox{-}}#1\operatorname{\mbox{-}trunc\mbox{-}map␣}}}

\newcommand{\tr}[3]{\ensuremath{\operatorname{tr}^{#1}_{#2}#3}}
\newcommand{\fib}{\ensuremath{\operatorname{fiber}}}
\newcommand{\refl}{\ensuremath{\operatorname{refl}}}
\newcommand{\id}{\ensuremath{\operatorname{id}}}
\newcommand{\idequiv}{\ensuremath{\operatorname{id\mbox{-}equiv}}}
\newcommand{\reflhtpy}{\ensuremath{\operatorname{refl\mbox{-}htpy}}}
\newcommand{\ap}[1]{\ensuremath{\operatorname{ap}_{#1}}}

\newcommand{\total}{\ensuremath{\operatorname{tot}}}

\newcommand{\image}{\ensuremath{\operatorname{image}}}
\newcommand{\incl}{\ensuremath{\operatorname{incl}}}

\newcommand{\ittype}[1]{\ensuremath{\textnormal{is-iter-}#1\textnormal{-type}}}

\newcommand{\T}{\ensuremath{\operatorname{P}}}
\newcommand{\V}{\ensuremath{\operatorname{V}}}
\newcommand{\W}[2]{\ensuremath{\operatorname{W}_{#1}#2}}
\newcommand{\M}[2]{\ensuremath{\rotatebox[origin=c]{-180}{$\operatorname{W}$}_{#1}#2}}
\newcommand{\desup}{\ensuremath{\operatorname{desup}}}
\newcommand{\ssup}{\ensuremath{\operatorname{sup}}}

\newcommand{\El}{\ensuremath{\operatorname{El}}}

\newcommand{\Mm}{\ensuremath{\mathrm{M}}}

\newcommand{\E}{\ensuremath{\operatorname{Unfold}}}
\newcommand{\F}{\ensuremath{\mathrm{F}}}

\newcommand{\source}{\ensuremath{\operatorname{source}}}
\newcommand{\target}{\ensuremath{\operatorname{target}}}
\newcommand{\Target}{\ensuremath{\operatorname{Target}}}

\newcommand{\Bisim}[2]{\ensuremath{#1\text{-}\mathrm{Bisim}_{#2}}}

\newcommand{\Coalg}[1]{\ensuremath{#1\text{-}\mathrm{Coalg}}}

\newcommand{\Hom}[1]{\ensuremath{\operatorname{Hom}_{#1}}}
\newcommand{\idecoalg}{\ensuremath{\operatorname{δ}}}

\newcommand{\longunderscore}{\rule{6pt}{0.4pt}}

\newcommand{\iscoittypen}[2]{\ensuremath{\textnormal{is-coit-}#1\textnormal{-type}_{#2}}}
\newcommand{\iscoittype}[1]{\ensuremath{\textnormal{is-coit-}#1\textnormal{-type}}}
\newcommand{\corec}{\ensuremath{\operatorname{corec}}}

\newcommand{\tar}{\ensuremath{\operatorname{t}}}

\definecolor{dkblue}{rgb}{0,0.1,0.5}
\definecolor{lightblue}{rgb}{0,0.5,0.5}
\definecolor{dkgreen}{rgb}{0,0.6,0}
\definecolor{dkbrown}{rgb}{0.4,0,0}
\definecolor{dkviolet}{rgb}{0.6,0,0.8}

\newcommand{\todo}[1]{}
\newcommand{\niccolo}[1]{}
\newcommand{\elisabeth}[1]{}
\newcommand{\hakon}[1]{}
\newcommand{\philipp}[1]{}


\begin{document}

\title[Terminal Coalgebras and Non-wellfounded Sets in HoTT]{Terminal
Coalgebras and Non-wellfounded Sets\texorpdfstring{\\}{} in Homotopy Type Theory}

\author[H. R.~Gylterud]{Håkon Robbestad
Gylterud\lmcsorcid{0009-0009-9871-1110}}[a]
\author[E.~Stenholm]{Elisabeth
Stenholm\lmcsorcid{0009-0001-5852-9594}}[a]
\author[N.~Veltri]{Niccolò Veltri\lmcsorcid{0000-0002-7230-3436}}[b]

\address{Department of Informatics, University of
Bergen, Thormøhlensgate 55B, 5006 Bergen, Norway}
\email{hakon.gylterud@uib.no, elisabeth@stenholm.one}
\address{Department of Software Science, Tallinn University of
Technology, Akadeemia tee 15a, 12618 Tallinn, Estonia}
\email{niccolo@cs.ioc.ee}

\begin{abstract}
  Non-wellfounded material sets have previously been modelled in
  Martin-Löf type theory by Lindström using setoids. In this paper we
  construct models of non-wellfounded material sets in Homotopy Type
  Theory (HoTT) where equality is interpreted as the identity type. The
  first model satisfies Scott's anti-foundation Axiom (SAFA) and
  dualises the construction of iterative sets. The second model
  satisfies Aczel's anti-foundation Axiom (AFA), and is constructed by
  adapting Aczel--Mendler's terminal coalgebra theorem to type theory,
  which requires propositional resizing.

  In a bid to extend coalgebraic theory and anti-foundation axioms to
  higher type levels, we formulate generalisations of AFA and SAFA, and
  construct a hierarchy of models which satisfy the SAFA
  generalisations. These generalisations build on the framework of
  Univalent Material Set Theory, previously developed by two of the
  authors.

  Since the model constructions are based on M-types, the paper also
  includes a characterisation of the identity type of M-types as indexed
  M-types.

  Our results are formalised in the proof-assistant Agda.
\end{abstract}

\maketitle

\section{Introduction}\label{introduction}

In non-wellfounded set theory, the concept of a material set is expanded
beyond the cumulative hierarchy. The allowance for non-wellfounded sets,
such as the Quine atom \(q≔\{q\}\), makes it easier to study circular
phenomena and structures such as transitions systems and
streams\footnote{Aczel~\cite[Chapter 8]{aczel1988} gives an introduction to
  applications of non-wellfounded sets.}. In what follows, we seek to
integrate non-wellfounded set theory into Homotopy Type Theory
(HoTT)---a relatively new framework for mathematics, which supports
higher dimensional structures as first-class citizens with the powerful
\emph{Univalence Axiom} and higher inductive types
\citep[Chapter 6]{hottbook}. Our aim is to take classical notions from
universal coalgebra and non-wellfounded set theory and extend them to
higher-dimensional structures.

Wellfounded material set theory has been studied in Martin-Löf type
theory since 1978 with the introduction of Aczel's setoid model of
Constructive Zermelo--Fraenkel set theory~(CZF)~\citep{aczel1978}.
Non-wellfounded set theory in Martin-Löf type theory was studied already
in 1989 by Lindström, when she constructed a setoid based model of
constructive ZF⁻ (ZF without the axiom of foundation) + Aczel's
anti-foundation axiom (AFA) \citep{lindstrom-1989}.

These two models of material set theory were, as mentioned, setoid
based, meaning that equality was interpreted as a binary relation
distinct from Martin-Löf's identity type. This was rectified in the
model presented in the HoTT Book \citep{hottbook}, which constructed a
model of wellfounded set theory using a higher inductive type, in which
equality was interpreted as the identity type.

Gylterud~\cite{gylterud-iterative} then constructed a model, \((V⁰,∈)\),
equivalent to the HoTT Book model, but which did not require higher
inductive types for the construction. This construction and its
properties have been further explored by Gratzer et al.~\cite{gratzer2024category}.
One important aspect of \(V⁰\) is its role as the initial algebra of the
\(U\)-restricted powerset functor \(\T⁰_U : \Type → \Type\), which maps
\(X ↦ ∑_{A:U} \left(A ↪ X\right)\). One of the ideas we explore here is
to construct the terminal coalgebra for \(\T⁰_U\) to use as a model of
non-wellfounded sets, filling out the question mark in the table below.

\begin{longtable}[]{@{}lll@{}}
\toprule\noalign{}
& Setoid & Identity type \\
\midrule\noalign{}
\endhead
\bottomrule\noalign{}
\endlastfoot
Foundation & Aczel 1978 & Gylterud 2018 \\
Anti-foundation & Lindström 1989 & ? \\
\end{longtable}

We show that the terminal coalgebra for \(\T⁰_U\) would indeed yield a
model of Aczel's anti-foundation axiom (AFA):

\begin{quote}
\textbf{AFA}: Any (directed) graph can be uniquely decorated with sets
such that elementhood between the sets coincides with edges in the
graph. \citep{aczel1988}
\end{quote}

As we shall see, the classical Aczel--Mendler construction
\citep{aczel-mendler} can be adapted to the HoTT setting and constructs
a terminal coalgebra for \(\T⁰_U\), but it requires propositional
resizing---an impredicative axiom \citep[p. 116]{hottbook}.

In addition to the Aczel--Mendler construction, we provide a new
construction, \(V⁰_∞\), of non-wellfounded sets in HoTT which dualises
the construction of \(V⁰\), but which surprisingly \emph{does not} yield
a terminal coalgebra for \(\T⁰_U\). It is a third fixed point---neither
initial nor terminal. This type is a model of Scott's anti-foundation
axiom (SAFA, described by~Aczel~\cite[p.45 and p.49]{aczel1988}\footnote{Aczel
  attributes this axiom to Scott~\cite{scott1960}, an unpublished preprint
  which the authors at the time of writing have not gotten a hold of.}),
an alternative anti-foundation axiom to AFA. SAFA is based on the
concept of \emph{Scott extensionality}. A graph is Scott extensional if
equality of nodes in the graph coincides with isomorphism of unfolding
trees.

\begin{quote}
\textbf{SAFA}: Every Scott extensional graph can be
injectively\footnote{A decoration is \emph{injective} if an equality
  between the sets decorating nodes implies equality of the nodes being
  decorated.} decorated with sets, and the graph of all sets with edges
symbolising elementhood is Scott extensional.
\end{quote}

\textbf{Remark:} This formulation of SAFA is one of the variations
considered in Aczel's book, where it goes by the name of
AFA\(^{{≅}^t}\). This is the special case of AFA\(^{\sim}\), defined on
page 45 of Aczel's book, for the relation \({{≅}^t}\), which Aczel
introduced on page 49.

We also explore possible extensions of anti-foundation axioms to higher
types. In HoTT, there is a fundamental notion of \(n\)-type arising from
the iterative application of identity
types~\citep[Definition 7.1.1]{hottbook}. The 0-types are the sets,
where much of classical mathematics takes place. But even for
down-to-earth mathematics such as combinatorics, higher types can play a
role. Groupoids, that is 1-types, show up for instance in Joyal's theory
of combinatorial species \citep{bergeron1998}. We therefore propose
generalisations, \(n\)-AFA and \(n\)-SAFA, of both AFA and SAFA to
\(n\)-types. The model construction \(V⁰_∞\) is presented as a general
construction, \(Vⁿ_∞\), which then satisfies \(k\)-SAFA for each
\(k≤n\).

The construction of \(Vⁿ_∞\) is based on M-types. These types were
constructed in HoTT by Ahrens et al.~\cite{ahrens-2015}. We provide some further
general results about M-types. In particular, we fully characterise the
identity types of M-types as indexed M-types.

\subsection{Related work}\label{related-work}

Iterated, non-wellfounded multisets were originally studied, in a
classical setting, by D’Agostino and Visser~\cite{finality-regained}. Similar to our
development, they fit Scott's anti-foundation axiom into a coalgebraic
setting, where sets are considered to be the unisets (i.e.~multisets
with, coiteratively, only one occurrence of each element) inside the
terminal coalgebra of an endofunctor---in their case on classes of sets.
Theorem \ref{thm:identity-on-m-types} in our work is closely related to
their result that Γ-bisimulation corresponds to isomorphic unfolding
trees. The functor Γ in their setting, is closely related to \(\T^∞\) in
our setting. Our results generalise this classical paper in a few
different directions. Firstly, we work in a constructive setting
(D’Agostino and Visser~\cite{finality-regained} uses law of the excluded middle in the proof
of Theorem 4.6 and the axiom of choice in Lemma 3.5). Secondly, we work
in HoTT so that our results can be interpreted at different homotopy
levels. And finally some of our results apply to polynomial functors in
general.

Terminal coalgebras of polynomial functors, a.k.a.~M-types, were
constructed in HoTT by Ahrens et al.~\cite{ahrens-2015}. Their iterative construction
cannot be straightforwardly used to build the terminal coalgebra of
other functors in HoTT, such as the finite powerset and the finite
multiset functors, as discussed by Veltri~\cite{veltri-2021} and
Joram and Veltri~\cite{joram-2023}. A method for constructing terminal coalgebras of a
large class of functors, subsuming the latter two complicated cases, has
been developed by Kristensen et al.~\cite{kristensen-2022}, building on preliminary
investigations by Møgelberg and Veltri~\cite{mogelberg-2019}, in an extension of Cubical
Type Theory with multi-clock guarded recursion enabling the
specification of coinductive types. Møgelberg and Veltri~\cite{mogelberg-2019} also
characterises the identity type of terminal coalgebras as the terminal
bisimulation, generalising the statement of our Theorem
\ref{thm:identity-on-m-types} to more general functors than polynomial
functors. It is however unclear whether these results can be replicated
internally in HoTT without the extension of the type system with
multi-clock guarded recursion.

In the setting of Pure Type Systems (PTS), there is another connection
between type theory and non-wellfounded set theory: the theory IZ + AFA
+ TC (Intuitionistic Zermelo with Aczel's antifoundation axiom and
transitive closures) is equiconsistent with a PTS called λZ
\citep{miquel2004}. This equiconsistency is established via an
interpretation using pointed graphs in λZ. Unlike Lindström's model
\citep{lindstrom-1989}, or our own, this interpretation is a syntactic
translation, rather than an internal type of sets.

\newpage

\subsection{Contributions}\label{contributions}

The main contributions of this paper are the following:

\begin{itemize}
\tightlist
\item
  Construction of a fixed point for each of the non-polynomial functors
  \(X ↦ ∑_{A : U} \left(A ↪ₙ X\right)\), which is distinct from both the
  initial algebra and the terminal coalgebra.
\item
  Adapting the Aczel--Mendler construction \citep{aczel-mendler} to type
  theory, assuming propositional resizing.
\item
  Applying the HoTT version of Aczel--Mendler to construct a terminal
  coalgebra for the \(U\)-restricted powerset functor.
\item
  A demonstration that this terminal coalgebra yields a model of set
  theory validating Aczel's anti-foundation axiom, with the identity
  type serving as equality.
\item
  Showing that Scott's anti-foundation axiom has a constructive model in
  HoTT, with the identity type as equality.
\item
  A characterisation of the identity types of M-types as indexed
  M-types.
\end{itemize}

\subsection{Formalisation}\label{formalisation}

The results in this paper have been formalised in the Agda proof
assistant~\citep{agda} and has been type checked with version 2.8.0. The
formalisation of Section \ref{coalgebra} to Section
\ref{the-coiterative-hierarchy} builds on the agda-unimath library
\citep{agda-unimath}, which is an extensive library of formalised
mathematics from the univalent point of view. The results in Section
\ref{the-terminal-tux2070-coalgebra} are formalised using Cubical
Agda---an extension of Agda with features from cubical type theory
\citep{cubical-agda}. This section has been type checked with Cubical
Agda version 0.9. Although the formalisation uses Cubical Agda, the
informal proofs in Section \ref{the-terminal-tux2070-coalgebra} are
carried out in the same framework as the rest of the article.

The formalisation of Sections
\ref{coalgebra}--\ref{the-coiterative-hierarchy} in this paper have been
included in a larger library on material set theory in HoTT, which can
be found here: \url{https://git.app.uib.no/hott/hott-set-theory}. As the
formalisation is structured slightly differently than the outline of
this paper, there are a few results which do not have an exact
counterpart in the code base. All these results are simple corollaries
or variations of results which have been formalised. Importantly, all
the main results are fully formalised. The formalisation of Section
\ref{the-terminal-tux2070-coalgebra} can be found at:
\url{https://github.com/niccoloveltri/aczel-mendler}. Throughout the
paper there will also be clickable links to specific lines of Agda code
corresponding to a given result. These will be shown as the Agda logo
\AgdaHen.

\subsection{Outline of the paper}\label{outline-of-the-paper}

To help the reader navigate, here is a short outline of the sections
following this introduction:

\textbf{Section \ref{coalgebra}} develops some basic notions of
coalgebra for wild endofunctors on \(\Type\). This provides a background
on working with coalgebra in Homotopy Type Theory and goes into some
detail on the identity type of terminal algebras for polynomial
endofunctors on \Type ~(i.e.~M-types). This gives a characterisation of
said identity types (Theorem \ref{thm:identity-on-m-types}), which is
used for our models in later chapters which are based on M-types.

\textbf{Section \ref{material-set-theory-in-homotopy-type-theory}}
recounts notions from \emph{Univalent Material Set Theory}
\citep{GylterudStenholm2026} which are used in this
article\hspace{1pt}---\hspace{1pt}in particular the notion of an
extensional ∈-structure and the coalgebraic perspective on it. We also
recall a generalisation of Rieger's theorem \citep{rieger1957}.

\textbf{Section \ref{afa-and-safa-in--structures}} provides the
generalised formulations of AFA and SAFA to ∈-structures in Homotopy
Type Theory. The generalisations are indexed by a truncation level,
where the \(0\)-AFA and \(0\)-SAFA are equivalent to the classical
formulations. The two coincide in the untruncated case: ∞-AFA is
equivalent to ∞-SAFA. In this section we also show, in Theorem
\ref{thm:AFA-from-terminal-coalg}, that terminal coalgebras of certain
generalised powerset functors, \(\Tⁿ\), are models of \(n\)-AFA.

\textbf{Section \ref{the-coiterative-hierarchy}} constructs a hierarchy
of models for non-wellfounded sets in Homomotopy Type Theory. These
model the generalised \(k\)-SAFA axioms (Theorem
\ref{thm:vninf-has-nsafa}). As types, these models are fixed points of
the generalised powerset functors, \(\Tⁿ\).

\textbf{Section \ref{the-terminal-tux2070-coalgebra}} describes a
general construction (Theorem \ref{thm:aczel-mendler}) of terminal
coalgebras for functors satisfying a certain accessibility condition, by
adapting the classical Aczel--Mendler
construction~\citep{aczel-mendler}. This establishes, in Corollary
\ref{cor:p0-has-terminal-coalgebra}, that the restricted powerset
functor has a terminal coalgebra in Homotopy Type Theory, assuming
propositional resizing. It remains an open question if such a terminal
coalgebra exists in Homotopy Type Theory without this assumption.

\subsection{Notation and conventions}\label{notation-and-conventions}

The notation throughout the paper will follow common practice in HoTT.
One important convention to note is that the type \(∑_{a : A} B␣a → X\)
should be read as \(∑_{a : A} \left(B␣a → X\right)\), rather than
\(\left(∑_{a : A} B␣a\right) → X\).

The ambient type theory is assumed to contain M-types. This is not a
very restrictive assumption as it has been shown by Ahrens et al.~\cite{ahrens-2015}
that M-types can be constructed from inductive types in HoTT.

Throughout the paper we will take the type of truncation levels to be
the type \(ℕ^∞_{-2}\), i.e.~the usual truncation levels, but with a
supremum, ∞, such that \(\|P\|_∞ ≡ P\). Moreover, for computations we
have \(∞ - 1 = ∞ = ∞ +
1\). We will also use \(ℕ^∞_{-1}\) for the subset of truncation levels
excluding \(-2\), and \(ℕ_{-2}\) and \(ℕ_{-1}\) for the ones further
excluding ∞.

We will also take liberties with coercions of subtypes into their
ambient type to enhance the readability of theorems and proofs. Since
the results are all formalised in Agda, we allow ourselves this
simplification without worry of any loss of rigour. The same goes for
using some essentially small types (i.e.~\(X : \Type\) for which there
is \(A:U\) such that \(A ≃ X\)) in some places instead of their small
replacements.

\section{Coalgebra}\label{coalgebra}

In this section we develop some basic notions of coalgebra theory with
focus on coalgebras for certain wild endofunctors on \(\Type\) and in
particular polynomial functors (defined by containers).

The notion of an \(\F\)-coalgebra is usually formulated for functors on
categories. In HoTT, there is a whole spectrum of notions of categories,
depending on if one wants univalence or would like to restrict the type
levels of homomorphism types or objects or both. Somewhere on this
spectrum we find the \emph{wild categories} \citep{kraus_path_2019},
where objects and homomorphisms can be of any type level and no
univalence is required. It turns out that in order for \(\F\)-coalgebras
to form a wild category, the endofunctor \(\F\) must satisfy some
additional coherences. These coherences are satisfied definitionally by
endofunctors on \(\Type\) induced by containers.

\subsection{\texorpdfstring{Coalgebras on
\(\Type\)}{Coalgebras on \textbackslash Type}}\label{coalgebras-on-type}

In this setting we will be interested in wild functors
\(\F : \Type → \Type\), which is an operation on types with an action
\((X → Y) → (\F␣X → \F␣Y)\), which we denote by juxtaposition \(\F␣f\),
which preserves composition and the identity function.

An \(\F\)-coalgebra is a pair \((A,α)\), where \(A : \Type\) and
\(α : A → \F␣A\). As is usual in universal coalgebra, we require no
comonadicity of \(\F\) nor coassociativity of \(α\) (i.e. α being an
algebra for \(\F\) as a comonad). We will also here settle on some
notation for standard notions of universal coalgebra, adapted to the
HoTT setting.

Given a wild functor \(\F\), coalgebras on \(\F\) form a wild category
if \(\F\) additionally satisfies some coherences. This generalises the
usual construction of the category of \(\F\)-coalgebras
\citep{Rutten2000} in an obvious way, and the coherences required fall
out when proving proof-relevant associativity. The proof below is a
sketch highlighting where the coherences are needed, the full details of
the proof can be found in the formalisation.

\begin{theorem}[The wild category of F-coalgebras \agdalink{https://elisabeth.stenholm.one/non-wellfounded-set-theory/v4/functor.coalgebra.html\#15092}]

\label{thm:f-coalgebra-cat}Let \(\F : \Type → \Type\) be a wild
endofunctor on the wild category of types and functions, satisfying the
following coherences:

\begin{enumerate}
\def\labelenumi{\arabic{enumi}.}
\tightlist
\item
  For any function \(f\), the two homotopies \(\F␣\id ∘ \F␣f ∼ \F␣f\)
  obtained from the fact that \(\F\) preserves identities and
  composition respectively, are homotopic.
\item
  For any function \(f\), the two homotopies \(\F␣f ∘ \F␣\id ∼ \F␣f\)
  obtained from the fact that \(\F\) preserves identities and
  composition respectively, are homotopic.
\item
  For any composable functions \(f\), \(g\) and \(h\), the two
  homotopies \(\F␣h ∘ \F␣g ∘ \F␣f ∼ \F␣(h ∘ g ∘ f)\) obtained from the
  fact that \(\F\) respects composition, associated in two ways, are
  homotopic.
\end{enumerate}

Then the type of \(\F\)-coalgebras
\[\Coalg{\F} ≔ ∑_{A : \Type} A → \F␣A\] and the type of \(\F\)-coalgebra
homomorphisms
\[\Hom{\Coalg{\F}}␣(A,α)␣(B,β) ≔ ∑_{f : A → B} β ∘ f ∼ \F␣f ∘ α\] form a
wild category.\end{theorem}

\begin{proof}

The underlying function of the identity homomorphism is the identity
function, and the homotopy is the proof that \(\F\) respects identities.
Composition of two homomorphisms is function composition of the
underlying functions and horizontal composition of the commuting
squares.

The identity laws hold definitionally on the underlying functions. The
non-trivial part of the laws are showing that the homotopies are equal.
Most of the steps follow from the groupoid laws of homotopies. For one
of the identity laws (depending on how composition of squares is
defined) one needs to use the fact that the two ways of composing
homotopies horizontally, differing in which homotopy is applied first,
are homotopic. Crucially, one needs the homotopy of 1. as a step in
showing that the identity is neutral with respect to postcomposition,
and, respectively, 2. for precomposition.

For showing that composition is associative, again, the law holds
definitionally for the underlying functions. For the homotopies, one
needs to use the same fact mentioned above about the two ways of
composing homotopies horizontally. Moreover, one needs the homotopy of
3. in one of the steps.\end{proof}

It is important to note that since the carrier of the codomain, \(B\),
can be of any type level, the second component of the type of
homomorphisms, namely \(β ∘ f
∼ \F␣f ∘ α\), is a structure, not just a property.

There are some definitions on coalgebras that will be useful throughout
the paper.

\begin{definition}

An \(\F\)-coalgebra \((A,α)\) is \textbf{extensional} if
\(α : A → \F␣A\) is an embedding.\end{definition}

Through the lens of type levels, we can also see a close connection
between two important properties of coalgebras, being \emph{terminal}
and being \emph{simple}:

\begin{definition}

Let \((A,α)\) be an \(\F\)-coalgebra.

\begin{itemize}
\tightlist
\item
  \((A,α)\) is \textbf{terminal} if for every \(\F\)-coalgebra,
  \((B,β)\), the type \(\Hom{\Coalg{\F}}␣(B,β)␣(A,α)\), is contractible.
\item
  \((A,α)\) is \textbf{simple} if for every \(\F\)-coalgebra, \((B,β)\),
  the type \(\Hom{\Coalg{\F}}␣(B,β)␣(A,α)\), is a proposition.
\end{itemize}

\end{definition}

\textbf{Remark:} In the classical literature \citep{Rutten2000} simple
is usually defined as: ``any outgoing epimorphism is an isomorphism''.
This is equivalent to any ingoing homomorphisms being equal, if the
functor preserves weak pullbacks. We choose here to use the second
characterisation as our definition, as this made the formulation of our
propositions more uniform. But the reader should be aware that they are
not equivalent without further assumptions.

The following is immediate from the definitions:

\begin{lemma}

A terminal \(\F\)-coalgebra is simple.\end{lemma}

\subsection{Bisimulation}\label{bisimulation}

Bisimulation is another central notion of coalgebra theory
\citep{Rutten2000}. In short, a bisimulation is just a span in the
category of \(\F\)-coalgebras, or a relation on the coalgebra that
relates elements in a way compatible with the coalgebra structure. More
specifically, a bisimulation on the \(\F\)-coalgebra \((X , m)\)
consists of an \(\F\)-coalgebra \((R , α)\) and functions
\(p₀, p₁ : R → X\) together with homotopies making the following diagram
commute:

\begin{center}
  \begin{tikzcd}
      X \arrow[dd, "m"] &  & R \arrow[ll, "p₀"' ] \arrow[rr, "p₁"] \arrow[dd, "α"] &  & X \arrow[dd, "m"] \\
                        &  &                                                       &  &                   \\
      \F␣X              &  & \F␣R \arrow[ll, "\F␣p₀"] \arrow[rr, "\F␣p₁"']         &  & \F␣X              
  \end{tikzcd}
\end{center}

A morphism of bisimulations is an F-coalgebra homomorphism
\(f : (R,α) → (R',α')\) between the underlying F-coalgebras, \((R,α)\)
and \((R',α')\), of the bisimulations, along with a filling of the left
and right triangular prisms of the resulting diagram:

\begin{center}
    \begin{tikzcd}
        & R \arrow[dd, "α" description, near start] \arrow[rr, "p₁" description] \arrow[ld, "p₀" description] \arrow[rd, "f" description] & & X \arrow[dd, "m"] \\
        X \arrow[dd, "m"'] & & R' \arrow[dd, "α'" description, near start] \arrow[ru, "p'₁" description] \arrow[ll, "p'₀" description] & \\
        & \F␣R \arrow[rr, "␣\F␣p₁"'] \arrow[ld, "\F␣p₀" description] \arrow[rd, "\F␣f" description]                        & & \F␣X \\
        \F␣X & & \F␣R' \arrow[ru, "\F␣p'₁" description] \arrow[ll, "\F␣p'₀"] &                         
    \end{tikzcd}
\end{center}

\begin{definition}[F-bisimulations on an F-coalgebra]

Let \(\F : \Type → \Type\) be a wild endofunctor and let \((X,m)\) be an
\(\F\)-coalgebra.

\begin{itemize}
\tightlist
\item
  \textnormal{(}\agdalink{https://elisabeth.stenholm.one/non-wellfounded-set-theory/v4/bisimulation.html\#1032}\textnormal{)}
  The type of \(\F\)-bisimulations on \((X,m)\) is the type
  \[\Bisim{\F}{(X,m)} ≔ ∑_{(R,α) : \Coalg{\F}} \Hom{\Coalg{\F}}␣(R,α)␣(X,m) × \Hom{\Coalg{\F}}␣(R,α)␣(X,m)\]
\item
  \textnormal{(}\agdalink{https://elisabeth.stenholm.one/non-wellfounded-set-theory/v4/bisimulation.html\#2371}\textnormal{)}
  Given two \(\F\)-bisimulations \(((R,α),σ₀,σ₁)\) and
  \(((R',α'),σ'₀,σ'₁)\), the type of \(\F\)-bisimulation homomorphisms
  from the first to the second is \begin{align*}
    \Hom{\Bisim{\F}{(X,m)}}&␣((R,α),σ₀,σ₁)␣((R',α'),σ'₀,σ'₁) \\
    &≔ ∑_{τ₀ : \Hom{\Coalg{\F}}␣(R,α)␣(R',α')}
        \left(σ'₀ ∘ τ₀ = σ₀\right) \\
        &\hspace{118pt}× \left(σ'₁ ∘ τ₀ = σ₁\right)
  \end{align*}
\end{itemize}

\end{definition}

When doing set level mathematics, a bisimulation homomorphism (the
homotopies being propositions) would simply be an \(\F\)-coalgebra
homomorphism from the first coalgebra to the second. But since we have
no restrictions on the type levels of the carrier types, we also need
coherences on the homotopies involved in the bisimulations.

In universal coalgebra, there are several equivalent\footnote{Assuming
  the functor preserves weak pullbacks.} formulations of being a simple
coalgebra~\citep{Rutten2000}. One of the equivalent formulations is that
the identity bisimulation is the terminal bisimulation. The definition
below is a strengthening of the classical definitions, allowing proof
relevant bisimulations and coalgebras with higher homotopies.

\begin{definition}[\agdalink{https://elisabeth.stenholm.one/non-wellfounded-set-theory/v4/bisimulation.html\#1360}]

\label{def:identity-bisim}Let \((X,m)\) be an \(\F\)-coalgebra. The
\textbf{identity bisimulation} on \((X,m)\) is simply \((X,m)\) together
with the identity homomorphism:
\[Δ_{(X,m)} ≔ ((X,m), \id_{(X,m)}, \id_{(X,m)}).\]\end{definition}

\begin{definition}

\label{def:simple-coalgebra}Let \((X,m)\) be an \(\F\)-coalgebra. We say
that \((X,m)\) is \textbf{bisimulation simple} if the identity
bisimulation is terminal. That is: for every bisimulation
\(((R,α),σ₀,σ₁)\) on \((X,m)\) the type
\[\Hom{\Bisim{\F}{(X,m)}}␣((R,α),σ₀,σ₁)␣Δ_{(X,m)}\] is
contractible.\end{definition}

We can relate equality of homomorphisms with bisimulation homomorphisms
into the identity bisimulation. This can then be used to show that being
simple and being bisimulation simple is equivalent.

\begin{lemma}[\agdalink{https://elisabeth.stenholm.one/non-wellfounded-set-theory/v4/bisimulation.html\#3004}]

\label{lem:hom-eq-bisim-hom}Assume that \(\F\) satisfies coherence 1 of
Theorem \ref{thm:f-coalgebra-cat}. Let \((X,m)\) be an \(\F\)-coalgebra
and let \(((R,α),σ₀,σ₁)\) be a bisimulation on \((X,m)\). There is an
equivalence: \[\Hom{\Bisim{\F}{(X,m)}}␣((R,α),σ₀,σ₁)␣Δ_{(X,m)}
≃ \left(σ₀ = σ₁\right).\]\end{lemma}

\begin{proof}

Let \(((R,α),σ₀,σ₁)\) be a bisimulation on \((X,m)\), where
\(σ₀≡(p₀,H₀)\) and \(σ₁≡(p₁,H₁)\). The type of bisimulation
homomorphisms from \(((R,α),σ₀,σ₁)\) to the identity bisimulation is the
type of fillings of the following diagram:

\begin{center}
  \begin{tikzcd}
      & R \arrow[dd, "α" description, near start] \arrow[rr, "p₁" description] \arrow[ld, "p₀" description] \arrow[rd, "f" description] & & X \arrow[dd, "m"] \\
      X \arrow[dd, "m"'] & & X \arrow[dd, "m" description, near start] \arrow[ru, "\id" description] \arrow[ll, "\id" description] & \\
      & \F␣R \arrow[rr, "␣\F␣p₁"'] \arrow[ld, "\F␣p₀" description] \arrow[rd, "\F␣f" description]                        & & \F␣X \\
      \F␣X & & \F␣X \arrow[ru, "\F␣\id" description] \arrow[ll, "\F␣\id"] &                         
  \end{tikzcd}
\end{center}

Using (one of) the identity laws on the category of \(\F\)-coalgebras we
can collapse the front two squares, and then use the fact that having a
homomorphism in the middle which is equal to \((p₀,H₀)\) and \((p₁,H₁)\)
is the same as the two being equal. (Note that we need here the fact
that \(\F\) satisfies coherence 1 of Theorem \ref{thm:f-coalgebra-cat}
for the needed identity law to hold.)

Or, presented in a different way, we have a chain of equivalences:
\begin{align*}
  &∑_{τ : \Hom{\Coalg{\F}}␣(R,α)␣(X,m)}
    \left(\id ∘␣τ = σ₀\right)
    × \left(\id ∘␣τ = σ₁\right)\\
  &≃ ∑_{τ : \Hom{\Coalg{\F}}␣(R,α)␣(X,m)}
      \left(τ = σ₀\right)
      × \left(τ = σ₁\right) \\
  &≃ \left(σ₀ = σ₁\right). \qedhere
\end{align*}\end{proof}

\begin{lemma}[\agdalink{https://elisabeth.stenholm.one/non-wellfounded-set-theory/v4/bisimulation.html\#5246}]

Assume that \(\F\) satisfies coherence 1 of Theorem
\ref{thm:f-coalgebra-cat}. An \(\F\)-coalgebra is bisimulation simple if
and only if it is simple.\end{lemma}

\begin{proof}

Let \((X,m)\) be an \(\F\)-coalgebra. We need to show that the type of
\(\F\)-coalgebra homomorphisms from \((Y,n)\) to \((X,m)\) is a
proposition for every \(\F\)-coalgebra \((Y,n)\), if and only if the
identity bisimulation \(Δ_{(X,m)}\) is terminal. By Lemma
\ref{lem:hom-eq-bisim-hom} we have, for any two \(\F\)-coalgebra
homomorphisms \((f,H)\) and \((g,K)\) from \((Y,n)\) to \((X,m)\), an
equivalence
\[\Hom{\Bisim{\F}{(X,m)}}␣((Y,n),(f,H),(g,K))␣Δ_{(X,m)} ≃ \left((f,H) = (g,K)\right).\]
The left hand side is contractible if and only if the right hand side is
contractible. Note that the right hand side being contractible for every
pair \((f,H)\) and \((g,K)\) means that the type of homomorphisms from
\((Y,n)\) to \((X,m)\) is a proposition. Thus \(Δ_{(X,m)}\) is the
terminal \(\F\)-bisimulation if and only if the type of \(\F\)-coalgebra
homomorphisms from \((Y,n)\) to \((X,m)\) is a proposition, i.e.~if and
only if \((X,m)\) is simple.\end{proof}

\begin{corollary}[\agdalink{https://elisabeth.stenholm.one/non-wellfounded-set-theory/v4/bisimulation.html\#7565}]

\label{cor:identity-is-terminal-bisim}Assume that \(\F\) satisfies
coherence 1 of Theorem \ref{thm:f-coalgebra-cat}. Let \((X,m)\) be a
terminal \(\F\)-coalgebra. Then \((X,m)\) is bisimulation simple,
i.e.~the identity bisimulation, \(Δ_{(X,m)}\), is the terminal
\(\F\)-bisimulation on \((X,m)\).\end{corollary}

\textbf{Remark:} This result can be seen as a generalisation and a
strengthening of Theorem 18 in Ahrens et al.~\cite{ahrens-2015}, the coinduction
proof principle. It is a generalisation as it also applies to
non-polynomial functors, and a strengthening as it shows that the
identity not only contains any other bisimulation, but that it is in
fact the terminal bisimulation. Later in the paper
(Theorem~\ref{thm:identity-on-m-types}) we will see that in the special
case of polynomial functors, the identity type is equivalently the
terminal coalgebra for a certain indexed functor.

\subsection{Characterisation of bisimulations of polynomial functors and
the identity type of
M-types}\label{characterisation-of-bisimulations-of-polynomial-functors-and-the-identity-type-of-m-types}

The M-types are a class of coinductive types, dual to the inductive
W-types. Intuitively, while the elements of W-types are wellfounded
trees with specified branching types, the M-types are the types of all
trees with that branching type. Formally, each M-type is the terminal
coalgebra of a polynomial functor which specifies the branching type. A
polynomial functor is one which is induced by a container
\citep{abbott_containers_2005,
altenkirch_indexed_2015}. Put simply, a polynomial functor
\(\Type→\Type\) is one of the form \(X ↦ ∑_{a:A} B␣a → X\), for some
\(A : \Type\) and \(B : A → \Type\). The data \(A\),\(B\) is called a
container and denoted \(A ⊲ B\). The functor \(X ↦ ∑_{a:A} B␣a → X\), as
induced by the container \(A ⊲ B\) is denoted by
\(⟦ A ◁ B ⟧ : \Type → \Type\). This wild endofunctor satisfies the
coherences required by Theorem \ref{thm:f-coalgebra-cat}. The M-type
\(\M{a:A}{B␣a} : \Type\) is the underlying type of the terminal
coalgebra of \(⟦ A ◁ B ⟧\) and its coalgebra map is denoted by:
\[\desup_{A,B} : \M{a:A}{B␣a} → ⟦ A ◁ B ⟧(\M{a:A}{B␣a}).\]

There are also indexed versions of polynomial functors, containers and
M-types. The indexed versions generalise from functors \(\Type → \Type\)
to functors \((I→\Type)→(J→\Type)\). An indexed polynomial functor maps
\(X ↦ λj.∑_{a:A_j} ∏_{b : B_ja}→X(w␣j␣b)\), for some \(A:J → \Type\) and
\(B: ∏_{j:J} A␣j → \Type\) and \(w : ∏_{j:J}∏_{a:Aj}B_ja→I\). The
functorial action sends a family of maps \(f : ∏_{i : I} X␣i → Y␣i\) to
the family \(λjλ(a,σ).(a,λb.f␣(w␣j␣b)␣(σ␣b))\). The data
\(A\),\(B\),\(w\) is called an \emph{indexed container}\footnote{Note
  that what we here call \emph{indexed container} is what
  Altenkirch et al.~\cite{altenkirch_indexed_2015} call a \emph{doubly indexed
  container}, which is \emph{not} the same as what they call indexed
  containers.} and is denoted by \(A⊲(B,w)\). The induced polynomial
functor is denoted by \(⟦ A ◁ (B,w) ⟧ : (I → \Type)→(J→\Type)\). The
indexed M-types are the terminal coalgebras for indexed polynomial
endofunctors, i.e. when \(I=J\).

Throughout the rest of this section, let \(A ◁ B\) be a container. For
convenience, we introduce some notation for \(⟦ A ◁ B ⟧\)-coalgebras.
This notation goes back to Aczel~\cite{aczel1978}, where it was applied to
its prototypical W-type, but we will use it for coalgebras in general.

\textbf{Notation:} Given \(m : X → ⟦ A ◁ B ⟧␣X\), and \(x : X\) we will
denote by \(\overline{x} : A\) and
\(\widetilde{x} : B␣\overline{x} → X\) the unique elements defined by
\(m␣x = (\overline{x}, \widetilde{x})\), that is
\(\overline{x} ≔ π₀␣(m␣x)\) and \(\widetilde{x} ≔ π₁␣(m␣x)\). This
notation suppresses the map \(m\), but it should be clear from the
context which map the notation refers to, whenever it is used. This
notation will also be used for \(\Type\) coalgebras
\(m : X → ∑_{I : \Type} I → X\).

The identity type of a W-type can be characterised inductively
\citep{gylterud-multisets}. For elements \(x, y : \W{a : A}{B␣a}\) there
is an equivalence:
\[\left(x = y\right) ≃ ∑_{p : \overline{x}=\overline{y}} ∏_{b : B␣\overline{x}}
\widetilde{x}␣b = \widetilde{y}␣(\tr{B}{p}{b}).\]

The goal of this subsection is to give a similar characterisation of the
identity type of M-types: The identity type between two elements of an
M-type is an indexed M-type (Theorem \ref{thm:identity-on-m-types}).
This characterisation is slightly more involved than the one for
W-types, which was proved by straightforward induction, and goes through
some results of bisimulation theory.

This result is not surprising, but is very useful for working with
M-types in HoTT. When we later construct a model of Scott's
non-wellfounded sets, this characterisation is critical in proving local
smallness of the model. Furthermore, the characterisation of the
identity type follows from a characterisation of bisimulations of
polynomial functors as coalgebras for a related indexed polynomial
functor.

\begin{definition}[\agdalink{https://elisabeth.stenholm.one/non-wellfounded-set-theory/v4/container.bisimulation.html\#5368}]

\label{def:E-functor}Given an \(⟦ A ◁ B ⟧\)-coalgebra \((X,m)\), we
define the \(\left(X × X\right)\)-indexed polynomial functor
\begin{align*}
     &\E_{(X,m)} : \left(X × X → \Type\right) → \left(X × X → \Type\right) \\
     &\E_{(X,m)}␣R␣(x, y) ≔ ∑_{p : \overline{x}=\overline{y}} 
         ∏_{b : B␣\overline{x}} R␣(\widetilde{x}␣b, \widetilde{y}␣(\tr{B}{p}{b})).
 \end{align*} The functorial action is postcomposition on the second
component. The functor respects identities and composition
definitionally.\end{definition}

The subscript \((X,m)\) will sometimes be omitted if it is clear from
the context.

We can think of \(\E_{(X,m)}\) as unfolding a relation one step as
though it was a bisimulation. A coalgebra for this functor, as we will
see, is thus a bisimulation on \((X,m)\).

\begin{definition}

Given an \(⟦ A ◁ B ⟧\)-coalgebra \((X,m)\), we define the following
types:

\begin{itemize}
\tightlist
\item
  The type of \(\E_{(X,m)}\)-coalgebras is
  \[\Coalg{\E_{(X,m)}} ≔ ∑_{R : X × X → \Type} ∏_{(x,y) : X × X} R␣(x,y) → \E␣R␣(x,y).\]
\item
  Given two \(\E_{(X,m)}\)-coalgebras \((R,α)\) and \((Q,β)\), the type
  of \(\E_{(X,m)}\)-coalgebra homomorphisms between them is
  \begin{align*}
    &\Hom{\Coalg{\E_{(X,m)}}}␣(R,α)␣(Q,β) \\
    &\hspace{20pt}≔ ∑_{f : ∏_{(x,y) : X × X} R␣(x,y) → Q␣(x,y)} 
        ∏_{(x,y) : X × X} β␣(x,y) ∘ f␣(x,y) ∼ \E␣f␣(x,y) ∘ α␣(x,y).
  \end{align*}
\end{itemize}

\end{definition}

The identity type is an \(\E_{(X,m)}\)-coalgebra, for any pair
\((X,m)\).

\begin{definition}[\agdalink{https://elisabeth.stenholm.one/non-wellfounded-set-theory/v4/container.bisimulation.html\#5566}]

\label{def:idecoalg}Define the following map by path induction:
\begin{align*}
    &\idecoalg : ∏_{(x,y) : X × X} x = y → \E_{(X,m)}␣({=})\;(x,y) \\
    &\idecoalg␣(x,x)␣\refl ≔ (\refl, \reflhtpy).
\end{align*} The pair \((=,\idecoalg)\) is the \textbf{identity
\(\E_{(X,m)}\)-coalgebra}.\end{definition}

An equivalence between two type families \((A,P)\) and \((B,Q)\) where
\(A,B : \Type\) and \(P : A × A → \Type\) and \(Q : B × B → \Type\), is
a pair \((α,σ)\) where \(α : A
≃ B\) and \(σ : ∏_{(a,a') : A × A} P␣(a,a') ≃ Q␣(α␣a,α␣a')\). By
univalence, we can transfer results about one family along such an
equivalence to a result about the other family. In our case, there is an
equivalence between \(\E\)-coalgebras and homomorphisms, and
bisimulations and homomorphisms for polynomial functors.

For this equivalence, we need to introduce notation for the
\textbf{total space} of a relation \(R : X × X → \Type\). This is the
type of all pairs that are related by \(R\):
\[|R| ≔ ∑_{(x,y) : X × X} R␣(x,y).\] Moreover, a fiberwise map
\(g : ∏_{(x,y) : X × X} R␣(x,y) → R'␣(x,y)\), induces a map on the total
spaces: \begin{align*}
    &\total␣g : |R| → |R'| \\
    &\total␣g␣((x,y), r) ≔ ((x,y), g␣(x,y)␣r).
\end{align*}

We are now ready to show the aforementioned equivalence.

\begin{theorem}[\agdalink{https://elisabeth.stenholm.one/non-wellfounded-set-theory/v4/container.bisimulation.html\#27558}]

\label{thm:E-coalg-bisim-equiv}Let \((X,m)\) be an
\(⟦ A ◁ B ⟧\)-coalgebra. There is an equivalence of type families
between \[\left(\Coalg{\E_{(X,m)}}, \Hom{\Coalg{\E_{(X,m)}}}\right)\]
and
\[\left(\Bisim{⟦ A ◁ B ⟧}{(X,m)}, \Hom{\Bisim{⟦ A ◁ B ⟧}{(X,m)}}\right).\]\end{theorem}

\begin{proof}

We start by constructing an equivalence
\[e : \Coalg{\E_{(X,m)}} ≃ \Bisim{⟦ A ◁ B ⟧}{(X,m)}.\] To this end,
first we note the following chain of equivalences: \begin{align*}
    \E_{(X,m)}␣R␣(x,y) &\hspace{12cm}
\end{align*} \vspace{-10mm} \begin{align*}
        &≃ ∑_{p : \overline{x} = \overline{y}} ∑_{φ₁ : B␣\overline{x} → X}
            (\widetilde{y} ∘ \tr{B}{p}{} = φ₁) × \left(∏_{b : B␣\overline{x}} R␣(\widetilde{x}␣b,φ₁␣b)\right) \\
        &≃ ∑_{φ₁ : B␣\overline{x} → X} (m␣y = (\overline{x},φ₁))
            × \left(∏_{b : B␣\overline{x}} R␣(\widetilde{x}␣b,φ₁␣b)\right) \\
        &≃ ∑_{a : A}\,∑_{φ₀ , φ₁ : B␣a → X} (m␣x = (a , φ₀)) × (m␣y = (a , φ₁)) ×
        \left(∏_{b : B␣a} R␣(φ₀␣b,φ₁␣b)\right) \\
        &≃ ∑_{(a , φ) : ⟦ A ◁ B ⟧␣|R|}
        (m␣x = ⟦ A ◁ B ⟧␣(π₀ ∘ π₀)␣(a , φ)) × (m␣y = ⟦ A ◁ B ⟧␣(π₁ ∘ π₀)␣(a , φ)).
\end{align*} Denote the equivalence above by \(e'\). Following an
element \((p,σ) : \E_{(X,m)}␣R␣(x,y)\) along \(e'\), we see that it is
mapped to
\[(\overline{x},λ␣b.((\widetilde{x}␣b,\widetilde{y}␣(\tr{B}{p}{b})),σ␣b)) : ⟦ A ◁ B ⟧␣|R|.\]

The equivalence above then gives us the desired equivalence \(e\):
\begin{align*}
    &\Coalg{\E_{(X,m)}} \\
        &\hspace{25pt}≃∑_{R : X×X → \Type}\, ∏_{(x,y) : X × X} R␣(x,y) \\
        &\hspace{50pt}→ ∑_{(a , φ) : ⟦ A ◁ B ⟧␣|R|}
        (m␣x = ⟦ A ◁ B ⟧␣(π₀ ∘ π₀)␣(a , φ)) \\
        &\hspace{108pt}× (m␣y = ⟦ A ◁ B ⟧␣(π₁ ∘ π₀)␣(a , φ)) \\
        &\hspace{25pt}≃∑_{R : X×X → \Type}\, ∑_{α : |R| → ⟦ A ◁ B ⟧␣|R|}
            \left(m ∘\,π₀ ∘ π₀ ∼ ⟦ A ◁ B ⟧␣(π₀ ∘ π₀) ∘ α\right) \\
            &\hspace{150pt}× \left(m ∘\,π₁ ∘ π₀ ∼ ⟦ A ◁ B ⟧␣(π₁ ∘ π₀) ∘ α\right) \\
        &\hspace{25pt}≃∑_{R : \Type}\, ∑_{p : R → X×X}\, ∑_{α : R → ⟦ A ◁ B ⟧␣R}
            \left(m ∘\,π₀ ∘ p ∼ ⟦ A ◁ B ⟧␣(π₀ ∘ p) ∘ α\right) \\
            &\hspace{153pt}× \left(m ∘\,π₁ ∘ p ∼ ⟦ A ◁ B ⟧␣(π₁ ∘ p) ∘ α\right) \\
        &\hspace{25pt}≃\Bisim{⟦ A ◁ B ⟧}{(X,m)}.
\end{align*} A pair \((R,f) : \Coalg{\E_{(X,m)}}\) is mapped to
\begin{align}
\label{coalg-bisim-equiv-comp}
&e␣(R,f) \nonumber \\
&\hspace{10pt}= \Big(\left(|R|, λ␣((x,y),r).(\overline{x},λ␣b.((\widetilde{x}␣b,\widetilde{y}␣(\tr{B}{π₀␣(f␣(x,y)␣r)}{b})),π₁␣(f␣(x,y)␣r)␣b))\right), \\
&\hspace{34pt}(π₀,\longunderscore), (π₁, \longunderscore)\Big), \nonumber
\end{align} where the homotopies have been left out for ease of
readability.

Now that we have an equivalence \(e\) on the base types, we need to
construct for any two pairs \((R,f), (R',f') : \Coalg{\E_{(X,m)}}\) an
equivalence
\[\Hom{\Coalg{\E_{(X,m)}}}␣(R,f)␣(R',f') ≃ \Hom{\Bisim{⟦ A ◁ B ⟧}{(X,m)}}␣(e␣(R,f))␣(e␣(R',f')).\]
To increase readability, let the following denote the components of
\(e␣(R,f)\) and \(e␣(R',f')\):

\begin{itemize}
\tightlist
\item
  \(α : |R| → ⟦ A ◁ B ⟧␣|R|\),
\item
  \(α' : |R'| → ⟦ A ◁ B ⟧␣|R'|\),
\item
  \(H₀ : m ∘\,π₀ ∘ π₀ ∼ ⟦ A ◁ B ⟧␣(π₀ ∘ π₀) ∘ α\),
\item
  \(H'₀ : m ∘\,π₀ ∘ π₀ ∼ ⟦ A ◁ B ⟧␣(π₀ ∘ π₀) ∘ α'\),
\item
  \(H₁ : m ∘\,π₁ ∘ π₀ ∼ ⟦ A ◁ B ⟧␣(π₁ ∘ π₀) ∘ α\)
\item
  \(H'₁ : m ∘\,π₁ ∘ π₀ ∼ ⟦ A ◁ B ⟧␣(π₁ ∘ π₀) ∘ α'\).
\end{itemize}

In other words,
\[e␣(R,f) ≡ ((|R|,α),(π₀,H₀),(π₁,H₁))\quad \text{and}\quad e␣(R',f') ≡ ((|R'|,α'),(π₀,H'₀),(π₁,H'₁)).\]
We have the following chain of equivalences: \begin{align*}
    &\Hom{\Coalg{\E_{(X,m)}}}␣(R,f)␣(R',f') \\
    &\hspace{25pt}≃ ∑_{g : ∏_{(x,y) : X × X} R␣(x,y) → R'␣(x,y)}
        \,∏_{(x,y) : X × X} \,∏_{r : R␣(x,y)} \\
        &\hspace{80pt}e'␣(f'␣(x,y)␣(g␣(x,y)␣r)) = e'␣(\E␣g␣(x,y)␣(f␣(x,y)␣r)) \\
    &\hspace{25pt}≃ ∑_{g : ∏_{(x,y) : X × X} R␣(x,y) → R'␣(x,y)}
        ∑_{K : α'∘\,\total␣g ∼ ⟦ A ◁ B ⟧␣(\total␣g)\,∘\,α} \\
            &\hspace{91pt}\left((H'₀ ∘ \total␣g) · (\ap{⟦ A ◁ B ⟧␣(π₀ ∘ π₀)} ∘␣K) = H₀ \right) \\
            &\hspace{80pt}× \left((H'₁ ∘ \total␣g) · (\ap{⟦ A ◁ B ⟧␣(π₁ ∘ π₀)} ∘␣K) = H₁ \right) \\
    &\hspace{25pt}≃ ∑_{g : ∏_{(x,y) : X × X} R␣(x,y) → R'␣(x,y)}
    ∑_{K : α'∘\,\total␣g ∼ ⟦ A ◁ B ⟧␣(\total␣g)\,∘\,α} \\
        &\hspace{91pt}\left((H'₀ ∘ \total␣g) · (\ap{⟦ A ◁ B ⟧␣(π₀ ∘ π₀)} ∘␣K) = \tr{λ␣h.m∘\,h∼⟦ A ◁ B ⟧␣h\,∘\,α}{\refl}{H₀} \right) \\
        &\hspace{80pt}× \left((H'₁ ∘ \total␣g) · (\ap{⟦ A ◁ B ⟧␣(π₁ ∘ π₀)} ∘␣K) = \tr{λ␣h.m∘\,h∼⟦ A ◁ B ⟧␣h\,∘\,α}{\refl}{H₁} \right) \\
    &\hspace{25pt}≃ ∑_{g : |R| → |R'|}\,∑_{p : π₀ ∘ π₀ ∘ g = π₀ ∘ π₀}\,∑_{q : π₁ ∘ π₀ ∘ g = π₁ ∘ π₀}
    ∑_{K : α'∘\,\total␣g ∼ ⟦ A ◁ B ⟧␣(\total␣g)\,∘\,α} \\
        &\hspace{91pt}\left((H'₀ ∘ \total␣g) · (\ap{⟦ A ◁ B ⟧␣(π₀ ∘ π₀)} ∘␣K) = \tr{λ␣h.m∘\,h∼⟦ A ◁ B ⟧␣h\,∘\,α}{p}{H₀} \right) \\
        &\hspace{80pt}× \left((H'₁ ∘ \total␣g) · (\ap{⟦ A ◁ B ⟧␣(π₁ ∘ π₀)} ∘␣K) = \tr{λ␣h.m∘\,h∼⟦ A ◁ B ⟧␣h\,∘\,α}{q}{H₁} \right) \\
    &\hspace{25pt}≃ \Hom{\Bisim{⟦ A ◁ B ⟧}{(X,m)}}␣(e␣(R,f))␣(e␣(R',f')).
\end{align*} As always, the full details can be found in the
formalisation.\end{proof}

Now we are ready to characterise the identity type on \(\M{a : A}{B␣a}\)
as an indexed \(\Mm\)-type. Recall the identity coalgebra
\(({=}, \idecoalg)\) given in Definition \ref{def:idecoalg}.

\begin{theorem}[\agdalink{https://elisabeth.stenholm.one/non-wellfounded-set-theory/v4/container.m-types.html\#3884}]

\label{thm:identity-on-m-types}The pair \(({=}, \idecoalg)\) is the
terminal
\(\E_{\left(\M{a : A}{B␣a}, \desup_{A,B}\right)}\)-coalgebra.\end{theorem}

\begin{proof}

By Theorem \ref{thm:E-coalg-bisim-equiv}, \((=,δ)\) is the terminal
\(\E_{\left(\M{a : A}{B␣a}, \desup_{A,B}\right)}\)-coalgebra if and only
if \(e␣(=,δ)\) is the terminal \(⟦ A ◁ B ⟧\)-coalgebra bisimulation on
\(\left(\M{a : A}{B␣a},\desup_{A,B}\right)\), where \(e\) is the
equivalence on the base types given in that theorem.

The value of \(e␣(=,δ)\) is given by (\ref{coalg-bisim-equiv-comp}). By
path induction, the coalgebra map in \(e␣(=,δ)\), is equal to the map
\begin{align*}
    &f : |=| → ⟦ A ◁ B ⟧␣|=| \\
    &f␣((x,x),\refl) ≔ (\overline{x},λ␣b.((\widetilde{x}␣b,\widetilde{x}␣b),\refl)).
\end{align*} Additionally, applying the equivalence
\(|=| ≃ \M{a : A}{B␣a}\), we get: \begin{align*}
    e␣(=,δ)
        &= \left(\left(|=|, f\right), (π₀,\longunderscore), (π₁, \longunderscore)\right) \\
        &= \left(\left(\M{a : A}{B␣a}, \desup_{A,B}\right), \id_{\left(\M{a : A}{B␣a}, \desup_{A,B}\right)}, \id_{\left(\M{a : A}{B␣a}, \desup_{A,B}\right)}\right) \\
        &≡ Δ_{\left(\M{a : A}{B␣a}, \desup_{A,B}\right)}
\end{align*} The full details of this computation can be found in the
formalisation.

By Corollary \ref{cor:identity-is-terminal-bisim}, the identity
bisimulation, \(Δ_{\left(\M{a : A}{B␣a},
\desup_{A,B}\right)}\), is the\linebreak terminal \(⟦ A ◁ B ⟧\)-bisimulation on
\(\left(\M{a : A}{B␣a}, \desup_{A,B}\right)\), and thus \((=,δ)\) is the
terminal\linebreak
\(\E_{\left(\M{a : A}{B␣a}, \desup_{A,B}\right)}\)-coalgebra.\end{proof}

\section{Material set theory in Homotopy Type
Theory}\label{material-set-theory-in-homotopy-type-theory}

In this section we revisit two perspectives on material set theory, both
generalised to higher homotopy levels: the notion of an ∈-structure and
the coalgebraic viewpoint. These two perspectives were developed in
previous work by two of the authors, titled \emph{Univalent Material Set
Theory} \citep{GylterudStenholm2026}. In univalent material set theory,
the elementhood relation, \(x ∈ y\), is not always a proposition, but
can be a type of any level. The elements of \(x ∈ y\) are considered
occurrences of \(x\) in \(y\). While from the coalgebraic perspective,
the same development can be seen as going from a subset perspective to
general fibrations by generalising the powerset functor.

\subsection{∈-structures and univalent material set
theory}\label{structures-and-univalent-material-set-theory}

The notion of an ∈-structure generalises the usual notion of an
extensional model of set theory to allow the underlying type to not just
be a set, but a type of higher level.

\begin{definition}[\agdalink{https://elisabeth.stenholm.one/non-wellfounded-set-theory/v4/e-structure.core.html\#1878}]

An \textbf{∈-structure} is a pair \((V,∈)\) where \(V : \Type\) and
\(∈␣: V → V → \Type\), which is \textbf{extensional}: for each
\(x,y : V\), the canonical map \(x = y → ∏_{z : V} {z ∈ x} ≃ {z ∈ y}\)
is an equivalence of types.\end{definition}

Extensionality ensures that the ∈-relation characterises equality of
sets (up to equivalence of types).

We can stratify ∈-structures based on the type level of the ∈-relation.

\begin{definition}[\agdalink{https://elisabeth.stenholm.one/non-wellfounded-set-theory/v4/e-structure.core.html\#4010}]

Given \(n : \Nat_{-2}\), an ∈-structure \((V,∈)\) is said to be of
\textbf{level (n+1)} if for every \(x,y:V\) the type \(x ∈ y\) is an
\(n\)-type.\end{definition}

The elements of a given set can also be collected to form a type.

\begin{definition}[\agdalink{https://elisabeth.stenholm.one/non-wellfounded-set-theory/v4/e-structure.core.html\#5019}]

Given an ∈-structure \((V,∈)\) we define the type family \begin{align*}
      &\El: V → \Type \\
      &\El␣a ≔ ∑_{x : V} x ∈ a.
 \end{align*}\end{definition}

Since \(V\) is a large type, \(\El␣a\) is a priori also a large type.
However, in many cases it is essentially small.

\begin{definition}[\agdalink{https://elisabeth.stenholm.one/non-wellfounded-set-theory/v4/e-structure.u-like.html\#1915}]

An ∈-structure \((V,∈)\) is \textbf{\(U\)-like} if the type \(\El␣a\) is
essentially \(U\)-small for all \(a : V\).\end{definition}

In this paper we will almost exclusively focus on the anti-foundation
axioms, but at times we will see some examples where we will use things
like the empty set, ∅, and paring/finite unordered tupling. In univalent
material set theory unordered tuples must be subscripted with their type
level. We will only use type level 0 and type level 1 in the examples,
so it is sufficient here to note that \(\{a₀,⋯,a_{n-1}\}_0\) is the
usual set theoretic tupling where repetition is ignored, while
\(\{a₀,⋯,a_{n-1}\}_1\) is multiset tupling where for instance
\(∅ ∈ \{∅,∅\}_1\) is a type with two elements. There is also the notion
of ordered pairing, but it is uniform in type level and consists of a
choice of embedding \(〈-,-〉 : V×V↪V\). See \emph{Univalent Material
Set Theory} \citep{GylterudStenholm2026} for details.

\begin{definition}[\agdalink{https://elisabeth.stenholm.one/non-wellfounded-set-theory/v4/e-structure.core.html\#7351}]

Given an ∈-structure, \((V,∈)\), \textbf{an ordered pairing structure}
on \((V,∈)\) is an embedding \(V × V ↪ V\).\end{definition}

We will use ordered pairs extensively when formulating anti-foundation
axioms. We will rely on the following proposition, which follows
immediately from the fact that ordered pairing is an embedding.

\begin{proposition}

\label{ordered-pair-prop}Being an ordered pair is a mere proposition:
for an ∈-structure \((V,∈)\) with ordered pairing structure \(〈-,-〉\),
the type \(∑_{a,b:V}〈a,b〉 = x\) is a proposition, for all
\(x : V\).\end{proposition}

\subsection{Coalgebraic view of set
theory}\label{coalgebraic-view-of-set-theory}

There is a coalgebraic viewpoint of material set theory, where one
replaces the usual \(∈\)-relation on \(V\) (classically the class of all
sets) with a coalgebra structure \(V → P(V)\) in the category of classes
and class functors. The functor \(P\) is the powerset functor on classes
which assigns to each class the class of subsets of the class. The axiom
of foundation says that \(V\) is the initial \(P\)-algebra, while
Aczel's anti-foundation axiom says that \(V\) is the terminal coalgebra.
Other \(P\)-coalgebras are known in set theory as \emph{set-like} models
of set theory, and the Mostowski collapsing theorem can be framed in
these terms. See for instance Paul Taylor's work on these topics
\citep{taylor2023}.

In \emph{Univalent Material Set Theory} \citep{GylterudStenholm2026},
two of the authors of the current paper developed this coalgebraic
viewpoint of material set theory inside HoTT, generalising it from sets
to types of arbitrary type levels. Since the models developed later use
this framework, we will quickly revisit the central definitions here.

The powerset functor on classes has a close correspondent in HoTT,
namely the \(U\)-restricted powerset functor:
\begin{align*}
    &\T⁰_U : \Type → \Type \\
    &\T⁰_U␣X ≔ ∑_{A : U} A ↪ X.
\end{align*}

The functorial action of \(\T⁰_U\) is taking the forward image along the
function: \[\T⁰_U f␣(A,v) = (\image (f∘v), \incl(f∘v)).\] By applying
the type theoretic replacement principle
\citep{rijke2017,rijke2019,rijke2022}, the image lands in \(U\) (and
thus the functorial action is well-defined) if the codomain of \(f\) is
locally \(U\)-small\footnote{A type A is \emph{locally U-small} if the
  identity type \(a=a'\) is essentially \(U\)-small for every
  \(a,a':A\).}. We will therefore restrict the application of this
functor to locally small types. By univalence, \(\T⁰_U\) preserves
local-smallness, hence one can regard it as a functor on locally small
types.

This notion of powerset is different from the one obtained by regarding
subtypes as maps into the type of \(U\)-small propositions. The two
notions coincide on types in \(U\), but differ on large types. In
particular, \(X ↦ \left(X →
\operatorname{hProp}_U\right)\) cannot have a fixed point, due to
Cantor's paradox. There is however no such obstacle for \(\T⁰_U\), which
is already known to have an initial algebra
\citep{hottbook,gylterud-iterative,GylterudStenholm2026}. As we shall
see later in this article, it also has a terminal coalgebra, assuming
propositional resizing, and a third fixed point (without assuming any
resizing). All fixed points are extensional coalgebras, which means that
they model the set theoretic extensionality axiom.

In univalent material set theory, one omits the requirement of having to
deal only with subtypes, and generalises to coalgebras for the
polynomial functor \(\T^∞_U\):
\begin{align*}
  &\T^∞_U : \Type → \Type \\
  &\T^∞_U␣X ≔ ∑_{A : U} A → X.
\end{align*}

The functorial action for \(\T^∞_U\) is simply postcomposition:
\[\T^∞_Uf(A,v)= (A,f∘v).\]

There is also a hierarchy of functors between \(\T⁰_U\) and \(\T^∞_U\),
where we restrict to \(n\)-truncated maps:
\begin{align*}
  &\T^{n+1}_U : \Type → \Type \\
  &\T^{n+1}_U␣X ≔ ∑_{A : U} A ↪_n X.
\end{align*}

The subscripted hooked arrow, \(A ↪_n X\), denotes an \(n\)-truncated
function \(A → X\). The \(n\) here ranges from \(-1\) to ∞, so that
\(\Tⁿ_U\) is defined for all \(n\) from \(0\) to ∞. The type
\(\T¹_U␣X\), for instance, is the type of set covers of \(X\).

The functorial action on \(\Tⁿ_U\) is taking \(n\)-images of the
composition: \[\Tⁿ_Uf(A,v) = (\image_n(f∘v),\incl_n(f∘v)).\] Just as for
\(\T⁰_U\), unless \(n=∞\), this is only well-defined on locally small
types.

Extensional coalgebras for these functors, that is, coalgebras for which
the coalgebra map is an embedding, correspond to ∈-structures in
univalent material set theory:

\begin{lemma}[∈-structures are coalgebras \agdalink{https://elisabeth.stenholm.one/non-wellfounded-set-theory/v4/e-structure.u-like.html\#8316}]
\label{coalgebra-structures}

For a fixed \(V\) and for \(n : \Nat^∞_{-1}\), having a \(U\)-like,
\((n+1)\)-level ∈-structure on \(V\) is equivalent to having a coalgebra
structure \(V ↪
 \T^{n+1}_U␣V\).\end{lemma}

\textbf{Remark:} This is Theorem 3 of Gylterud and Stenholm~\cite{GylterudStenholm2026}.   

\textbf{Notation:} As we do not work with several universes in this
article, we will often suppress mention of \(U\) in \(\Tⁿ_U\) and simply
write \(\Tⁿ\).

Since we will use it already in the definition of the anti-foundation
axioms, we will now take the opportunity to introduce the terminal
coalgebra of \(\T^∞\) which we will call \(\V^∞_∞\):
\[\V^∞_∞ ≔ \M{A:U}A.\]

This M-type comes equipped with a coalgebra structure
\(\desup^∞ : \V^∞_∞ →
\T^∞\V^∞_∞\), which is an equivalence, in a result analogous to Lambek's
lemma \citep{lambek1968}. Let \(\sup^∞ : \T^∞\V^∞_∞ → \V^∞_∞\) denote
the inverse of \(\desup^∞\). For any other \(\T^∞\)-coalgebra, \((X,m)\)
there is a unique coalgebra homomorphism
\(\corec^∞␣(X,m) : (X,m) → (\V^∞_∞ ,\desup^∞ )\). We will sometimes
suppress the coalgebra \((X,m)\) and only write \(\corec^∞\), when the
coalgebra is clear from the context.

\subsection{Fixed-point models}\label{fixed-point-models}

We have seen that ∈-structures are equivalent to extensional coalgebras
\(V ↪ \T^n\), (with \(n = ∞\) being the general case), but what if this
embedding is actually an equivalence? Then it turns out, in analogy with
Rieger's theorem \citep[Theorem III]{rieger1957} in classical set
theory, that generalisations of many constructive set theory axioms hold
in the corresponding ∈-structure, \((V,{∈})\). In particular, the
following axioms are shown to hold for fixed points, in
Gylterud and Stenholm~\cite{GylterudStenholm2026}:

\begin{itemize}
     \item Empty set.
     \item $U$-restricted $n$-separation.
     \item If $V$ is $(n+1)$-locally $U$-small, it has ∞-unordered $I$-tupling
          for all $(n-1)$-truncated types $I : U$.
     \item If $V$ is $(k+1)$-locally $U$-small, for some $k ≤ n$ then it has:
          \begin{itemize}
               \item $k$-unordered $I$-tupling for all $I : U$,
               \item $k$-replacement,
               \item $k$-union.
          \end{itemize}
     \item $V$ has exponentiation for all ordered pairing structures.
     \item $V$ has natural numbers represented by $f$ for any $(n-1)$-truncated
          representation $f : \Nat → V$.
\end{itemize}

The type level \(k\) on the axioms generalise from \(k=0\) of classical
set theory. In the case of a level 0 ∈-structure, one where \(x∈y\) is a
proposition for all \(x\) and \(y\), this specialises to the following
familiar axioms:

\begin{itemize}
     \item Empty set.
     \item $U$-restricted separation.
     \item If $V$ is locally $U$-small then it has:
          \begin{itemize}
               \item unordered pairs,
               \item replacement, and
               \item union.
          \end{itemize}
     \item $V$ has exponentiation for all ordered pairing structures.
     \item $V$ has natural numbers represented by $f$ for any of the usual representations.
\end{itemize}

Since the models we construct in this paper are all fixed points, we get
these basic axioms for free, provided we can prove that our models are
locally \(U\)-small. Proving this is the motivation for the previous
section on M-types. Thus, our focus will be on proving the axioms
particular to non-wellfounded set theory. But first, we will generalise
their formulation to the same level of generality as the other axioms of
univalent material set theory.

\section{AFA and SAFA in
∈-structures}\label{afa-and-safa-in--structures}

Most axioms of set theory, such as paring, union, separation and even
infinity, replacement or powerset, are \emph{set existence axioms} ---
they inform us which sets we can construct within the theory. All the
sets we can construct from these axioms alone are \emph{wellfounded}.
Classically, wellfounded sets are those without an infinite membership
chain: \[a₀ ∋ a₁ ∋ a₂ ∋ ⋯ \] Constructively, wellfoundedness is instead
formulated as an induction principle for ∈ or using an accessibility
predicate. In both constructive and classical traditions, the most
prominent theories include an axiom which states that, in fact, all sets
are wellfounded. This axiom is called regularity or \emph{foundation}.
It's a standard, classical result that the axiom of foundation is
independent of the rest \citep{bernays1954}. What is more, under certain
assumptions\footnote{The Axiom of Choice is more than sufficient, but
  the much milder axiom of wellfounded materialisation is enough
  \citep[cf. discussion in][after Lemma 6.46]{shulman_stack_2010}.} any
structure defined by sets can be defined by wellfounded sets.

When one removes the requirement that every material set must be
wellfounded, two questions arise:

\begin{enumerate}
\def\labelenumi{\arabic{enumi}.}
\tightlist
\item
  Which non-wellfounded sets exist?
\item
  When are two non-wellfounded sets equal?
\end{enumerate}

Anti-foundation axioms are properties of ∈-structures which give answers
to these two questions. In this text we consider two such axioms. The
first is Aczel's anti-foundation axiom (AFA), and the second is Scott's
anti-foundation axiom (SAFA). These answer the question slightly
differently, and in this section we will try to capture a formulation of
these in a way which generalises to ∈-structures of higher type levels.

The second question arises because extensionality does not fully
determine the equality between non-wellfounded sets. For instance, if
two sets satisfy the equations \(x = \{x,y\}₀\) and \(y = \{x\}₀\)
(using the notation introduced in Section
\ref{structures-and-univalent-material-set-theory}), both \(x=y\) and
\(x≠y\) are possible -- of course not in the same ∈-structure. The
0-subscript on the pairing is crucial, because if we used multiset
pairing, and let \(x = \{x,y\}₁\), it follows that \(x≠y\), since a pair
is never a singleton. This foreshadows the main thesis of this section,
that the difference between Aczel's and Scott's conceptions of
non-wellfounded sets is a matter of truncation level, from the
perspective of HoTT.

In elementary terms, AFA states that given any graph there is a unique
assignment of sets to the nodes of the graph, such that the elementhood
relations between the assigned sets coincides with the edges of the
graph. This gives both a way of constructing non-wellfounded sets (by
giving a graph) and a way of proving equalities between non-wellfounded
sets (showing that they can decorate the same node in a graph).

SAFA states that every graph where nodes have unique unfolding trees can
be decorated with sets (in the same sense as in AFA) and that for sets,
isomorphism of unfolding trees determines equality. Additionally, the
decoration is injective (since equality of nodes is determined by their
unfolding trees) and is unique among such decorations. This may at the
moment sound baroque and even ad hoc, but we will attempt to shed light
on this in Section \ref{subsec:safa}.

Why all these graphs? An answer to this question comes from universal
coalgebra. An ∈-structure being, in general a coalgebra for the functor
\(\T^∞\), and specifically a \(\Tⁿ\)-coalgebra in the case of
\(n\)-level structures \citep[cf.][Theorem 3]{GylterudStenholm2026}, the
non-well founded sets come from coalgebra maps into the structures. In
set-level mathematics, a graph is exactly a coalgebra \(X → \T⁰␣X\).
This emphasises looking at the out-edges from a node, and a coalgebra
map into an ∈-structure translates out-edges to elements. So, what we
will call a decoration of a graph is precisely a coalgebra homomorphism
from the induced coalgebra of the graph into the ∈-structure the graph
lives in.

\subsection{Graphs and decorations}\label{graphs-and-decorations}

Usually in mathematics, we think of graphs as structures consisting of
nodes and edges. However, in the formulation of the anti-foundation
axioms we will work with a slightly different notion of graph, as simply
a set of pairs. This leaves the domain of nodes implicit, which
simplifies the definition of a decoration. Another way of thinking of it
is that the domain of nodes in \(g\) is always the entirety of \(V\).

\begin{definition}[\agdalink{https://elisabeth.stenholm.one/non-wellfounded-set-theory/v4/e-structure.graphs.html\#1238}]

\label{def:graph}In an ∈-structure \((V,∈)\) with ordered pairing
structure \(〈-,-〉\), an element \(g : V\) is a \textbf{graph} if all
its elements are pairs. That is, there is a map
\[∏_{e : V} e ∈ g → ∑_{(x,y) : V × V} e = 〈x,y〉,\] or equivalently,
for every \(e : V\) such that \(e ∈ g\) there are \(\source␣e : V\) and
\(\target␣e : V\) such that
\(e = 〈\source␣e,\target␣e〉\).\end{definition}

\textbf{Remark:} The notation ``\(\source␣e\)'' and ``\(\target␣e\)''
suppresses mention of the specific proof element of \(e ∈ g\) which is
used to construct \(\source␣e\) and \(\target␣e\). However, this is
justified since ordered pairing is an embedding
\citep[Definition 5]{GylterudStenholm2026}, and hence
\(∑_{(x,y) : V × V} e = 〈x,y〉\) is a proposition. Thus any choice of
such a proof object yields equal results.

For the rest of the section we will assume that the ∈-structure we are
working with has an ordered pairing structure \(〈-,-〉\), in order to
avoid lengthy lists of assumptions in the statements of the results.

\begin{definition}[\agdalink{https://elisabeth.stenholm.one/non-wellfounded-set-theory/v4/e-structure.graphs.html\#1750}]

\label{def:target}Given a graph \(g : V\) in an ∈-structure \((V,∈)\),
define the type \(\Target␣g\), the subtype of \(V\) consisting of
targets of edges in \(g\), by
\(\Target␣g ≔ ∑_{y:V}∃_{x:V}〈x,y〉∈g\).\end{definition}

Since the domain of nodes in the graph is left implicit, a decoration
will be a universally defined function \(d : V →V\), where the
convention is that \(d␣x\) is empty if there are no edges
\(〈x,y〉 ∈ g\). When there is an edge \(〈x,y〉\) this edge should give
rise to an elementhood relation \(d␣y ∈ d␣x\). In fact, there should for
every \(z:V\) be an equivalence between \(z ∈ d␣x\) and the edges in
\(〈x,y〉 ∈ g\) for which \(z = d␣y\):

\begin{definition}[\agdalink{https://elisabeth.stenholm.one/non-wellfounded-set-theory/v4/e-structure.graphs.html\#2317}]

\label{def:n-decoration}For \(n : ℕ^∞_{-1}\), an
\textbf{\((n+1)\)-decoration} of a graph \(g : V\) in an ∈-structure
\((V,∈)\), is a map \(d
: V → V\) together with an element of the type
\[∏_{x,z : V} z ∈ d␣x ≃ \Bigg\|∑_{y : V} 〈x,y〉 ∈ g × d␣y = z\Bigg\|_{n}.\]\end{definition}

The truncation level restricts the level of \(d␣x\), so that, for
instance, in 0-level ∈-structures \(d␣x\) will be a set. The notion of
0-decoration is equivalent to the classical notion of decoration as a
function satisfying the equation
\(d(x) = \left\{\,d(y)\mid 〈x,y〉∈g\,\right\}₀\) \citep[cf.][Chapter
1]{aczel1988}. And, in terms of univalent material set
theory\footnote{See Gylterud and Stenholm~\cite{GylterudStenholm2026}, Definitions 7 and 8,
  for a discussion on n-truncated set comprehension and replacement.},
an \(n\)-decoration is a function satisfying the equation
\(d(x) = \left\{\,d(y)\mid 〈x,y〉∈ g\right\}_n\).

The notion of ∞-decoration is one where there is no truncation yielding
simply:
\begin{align}
\label{eq:infinity-decoration}
z ∈ d␣x &≃ \left(∑_{y : V} 〈x,y〉 ∈ g × d␣y = z\right).
\end{align}

Intuitively it says that \(d␣y\) occurs in \(d␣x\) precisely as many
times as \(〈x,y〉\) occurs in \(g\) (and that all elements of \(d␣x\)
are of the form \(d␣y\)).

There are two simple observations we can make if we know the level of
the ∈-structure.

\begin{itemize}
\tightlist
\item
  In an \(n\)-level ∈-structure, an \((n+1)\)-decoration is also an
  ∞-decoration since \(∑_{y : V} 〈x,y〉 ∈ g × d␣y = z\) has type level
  \(n\).
\item
  In an \(n\)-level ∈-structure, an ∞-decoration is also
  \(n\)-decoration, but the opposite is not always the case. For
  instance, in level \(0\), if \(d : V → V\) is an ∞-decoration, we know
  that \(∑_{y : V} 〈x,y〉 ∈ g × d␣y = z\) is a proposition since it is
  equivalent to \(z∈d␣x\) which is a proposition. Hence, the
  propositional truncation in the requirement for a \(0\)-truncation is
  superfluous and \(d\) is also a \(0\)-decoration. However, the graph
  \(g = \{〈a,b〉, 〈a,c〉\}₀\) cannot have an ∞-decoration in any
  \(0\)-level structure, if \(a\), \(b\) and \(c\) are distinct, since
  \(d␣b = d␣c = ∅\) and thus
  \(∅ ∈ d␣a ≃ \left(∑_{y : V} 〈a,y〉 ∈ g × d␣y = ∅\right) ≃ 2\), which
  is not a proposition. But, being wellfounded, \(g\) has a
  \(0\)-decoration, namely the one which assigns
  \(d␣x = \left\{\,∅\mid x = a\,\right\}₀\).
\end{itemize}

Classically, Scott's axiom is formulated in terms of injective
decorations, but we will instead use ∞-decorations as this generalises
to higher type levels. At level 0, the ∞-decorations are the injective
0-decorations. Note that, ``injective decoration'' does not mean that
\(d\) is injective on all of \(V\): Since \(g\) is a small set, \(d␣z\)
is ∅ on sets which are not nodes of \(g\) (i.e.~occurs in an edge in
\(g\)). But rather, what is meant by injective decoration is that it
becomes injective when restricted to the sets which are nodes in the
graph.

\subsection{\texorpdfstring{Coalgebraic characterisation of
\(n\)-decorations}{Coalgebraic characterisation of n-decorations}}\label{coalgebraic-characterisation-of-n-decorations}

Having seen that ∈-structures are the same as \(\T^∞\)-coalgebras, we
will now see that decorations can be identified with certain coalgebra
homomorphisms into these coalgebras. This is essentially what is proved
in Proposition \ref{prop:equiv-decoration-Tn-coalgebra} below and
mirrors the classical characterisation of decorations as coalgebra maps
into \(V\). However, to make the characterisation work, either the
functorial action must be adjusted for each \(n\), or the underlying
structure must be of level \(n\) (in the classical case \(n=0\)). We opt
to adjust the functorial action.

\begin{definition}

\label{def:Tinfn} Let \(n : ℕ^∞_{-1}\), and define a wild functor
\(\T^∞_{n+1} : \Type → \Type\) on types by
\(\T^∞_{n+1}␣X := ∑_{A:U} A→X\) and on functions by
\(\T^∞_{n+1}␣f␣(A,v) := (\image_n(f∘v),\incl_n(f∘v))\).\end{definition}

\textbf{Remark:} Notice that \(\T^∞_n\) is like a hybrid between
\(\T^∞\) and \(\T^n\): Since \(\T^∞_n\) and \(\T^∞\) have the same
action on types, a coalgebra for one is automatically a coalgebra for
the other. On the other hand, if two \(\T^∞_n\)-coalgebras factor
through \(\Tⁿ\)-coalgebras, the type of \(\T^∞_n\)-coalgebra
homomorphisms is equivalent to the type of \(\T^n\)-coalgebra
homomorphisms. The following commutative diagram summarises the
relationship between \(\T^n\) and \(\T^∞_n\). The unnamed arrows are the
\((n-1)\)-image map and the inclusion of \((n-1)\)-truncated functions
into functions.

\begin{center}
    \begin{tikzcd}
        \T^∞_n␣X \arrow[d,equal] \arrow[rrr,"\T^∞_nf"] & & & \T^∞_n␣Y \arrow[d,equal]\\
        \T^∞␣X \arrow[r,two heads] & \Tⁿ␣X \arrow[r,"\Tⁿf"] & \Tⁿ␣Y \arrow[r,hook] & \T^∞␣Y
    \end{tikzcd}
\end{center}

Let us for the rest of the subsection fix \(n : ℕ^∞_{-1}\) and a
\(U\)-like ∈-structure \((V,∈)\) and its associated \(\T^∞\)-coalgebra
structure \(m_{∈} : V →
\T^∞V\) (Lemma \ref{coalgebra-structures}). Assume also that \(V\) is
locally small and let \(x ≈ y\) denote the small type equivalent to the
identity type for \(x,y :
V\).

If we have a graph in \(V\), there are several ways of constructing a
coalgebra from it. Below, we define two closely related
\(\T^∞\)-coalgebra structures: \(m_g :V →\T^∞V\)~and
\(\tar_g : \Target␣g → \T^∞\left(\Target␣g\right)\), which will help
characterise decorations and define Scott's anti-foundation axiom.

\begin{proposition}[\agdalink{https://elisabeth.stenholm.one/non-wellfounded-set-theory/v4/e-structure.graph.to-P-n-coalgebra.html\#6530}]

\label{prop:graph-coalg} For each graph \(g : V\), there is a
\(\T^∞\)-coalgebra structure on \(V\) which we will call
\(m_g : V → \T^∞ V\) such that \(π₀(m_g␣x) ≃ ∑_{y:V}〈x,y〉∈g\) and
\(π₁(m_g␣x) : π₀(m_g␣x) → V\) becomes
\(π₀ : \left(∑_{y:V}〈x,y〉∈g\right) → V\) when transported along this
equivalence.\end{proposition}

\begin{proof}

Given \(x:V\) let
\(m_g␣x := (∑_{e : \overline{g}} \source␣(\widetilde{g}␣e) ≈ x,\target ∘\,\widetilde{g} ∘ π₀)\),
and observe that: \begin{align*}
    ∑_{e : \overline{g}} \source␣(\widetilde{g}␣e) ≈ x
        &≃ ∑_{y : V} ∑_{e : \overline{g}} 
            (\source␣(\widetilde{g}␣e) = x) × (\target␣(\widetilde{g}␣e) = y) \\
        &≃ ∑_{y : V} ∑_{e : \overline{g}} 
            〈\source␣(\widetilde{g}␣e),\target␣(\widetilde{g}␣e)〉 = 〈x,y〉 \\
        &≃ ∑_{y : V} \fib␣\widetilde{g}␣〈x,y〉 \\
        &≡ ∑_{y:V}〈x,y〉∈g.
\end{align*}

Note that the diagram

\begin{center}
    \begin{tikzcd}
        ∑_{e : \overline{g}} \source␣(\widetilde{g}␣e) ≈ x \arrow[rr, "≃"] \arrow[rd, "{\target ∘\,\widetilde{g} ∘ π₀}"'] &       & ∑_{y:V}〈x,y〉∈g \arrow[ld, "π₀"] \\
                                                   & V &                                                                                                 
    \end{tikzcd}
\end{center}

commutes up to definitional equality.\end{proof}

\textbf{Remark:} Ignoring size issues, justified by Proposition
\ref{prop:graph-coalg}, we will simply write:
\[m_g␣x = \left(∑_{y:V}〈x,y〉∈g,π₀\right).\] This is clearer to read
than coercing along an equivalence. A more careful treatment, without
notational abuse, is found in the formalisation.

\begin{lemma}[\agdalink{https://elisabeth.stenholm.one/non-wellfounded-set-theory/v4/e-structure.graph.to-P-n-coalgebra.html\#6425}]

\label{lma:mg-factor}If a graph \(g : V\) is an \(n\)-type in \((V,∈)\)
(i.e.~\(e∈g\) is an \((n-1)\)-type) then \(π₁␣(m_g␣x) : π₀␣(m_g␣x) → V\)
is \((n-1)\)-truncated, and hence \(m_g\) factors through a
\(\Tⁿ\)-coalgebra \(m_{n,g} : V → \TⁿV\).\end{lemma}

\begin{proof}

The map \(\target ∘\,\widetilde{g} ∘ π₀\) is \((n-1)\)-truncated since,
for any \(y : V\), we have the equivalences \begin{align*}
    \fib␣(\target ∘\,\widetilde{g} ∘ π₀)␣y
        ≃ \fib␣π₀␣y
        ≃ 〈x,y〉 ∈ g,
\end{align*} and the last type is \((n-1)\)-truncated. (The first
equivalence uses the commuting diagram in the proof of Proposition
\ref{prop:graph-coalg}.)\end{proof}

\begin{proposition}[\agdalink{https://elisabeth.stenholm.one/non-wellfounded-set-theory/v4/e-structure.graph.to-P-n-coalgebra.html\#8793}]

\label{prop:equiv-decoration-Tn-coalgebra} For each graph \(g:V\) there
is an equivalence between the type of \(n\)-decorations of \(g\) and the
type of \(\T^∞_n\)-coalgebra homomorphisms from \(m_g\) to
\(m_{∈}\).\end{proposition}

\begin{proof}

Given a graph \(g : V\) and a map \(d : V → V\) we have the following
chain of equivalences: \begin{align*}
    (m_{∈} ∘␣d ∼ \T^∞_n␣d ∘ m_g)
        &≃ ∏_{x : V} ∏_{z : V} \fib␣\widetilde{(d␣x)}␣z ≃ \fib␣(\incl_{n-1}␣(d ∘ \target ∘\,\widetilde{g} ∘ π₀))␣z \\
        &≃ ∏_{x : V} ∏_{z : V} z ∈ d␣x ≃ \Bigg\| ∑_{(y,p) : \fib␣d␣z} 〈x,y〉 ∈ g \Bigg\|_{n-1} \\
        &≃ ∏_{x : V} ∏_{z : V} z ∈ d␣x ≃ \Bigg\|∑_{y : V} 〈x,y〉 ∈ g × d␣y = z\Bigg\|_{n-1}
\end{align*} The first step uses Proposition 18 of
Gylterud and Stenholm~\cite{GylterudStenholm2026} which states that equality on slices is
equivalence on fibers. The second step uses the fact that the fiber of
the \((n-1)\)-image inclusion is the \((n-1)\)-truncation of the fiber
of the original map. The full equivalence then follows from the fact
that dependent sums preserve equivalences.\end{proof}

\begin{proposition}[\agdalink{https://elisabeth.stenholm.one/non-wellfounded-set-theory/v4/e-structure.graph.to-P-n-coalgebra.html\#7803}]

\label{prop:graph-coalg-Target} For each graph \(g : V\), the coalgebra
\(m_g\) restricts to \(\Target␣g\). We will call this coalgebra
structure \(\tar_g :
\Target g → \T^∞ \left(\Target␣g\right)\) and the subtype inclusion
\(π₀ : \Target g → V\) is a \(\T^∞\)-coalgebra
homomorphism.\end{proposition}

\begin{proof}

First, note that for any \(e : \overline{g}\),
\(\target␣(\widetilde{g}␣e)\) lies in \(\Target␣g\) as it is the child
of \(\source␣(\widetilde{g}␣e)\). Thus let
\(\tar_g(x,\longunderscore) = (∑_{e : \overline{g}} \source␣(\widetilde{g}␣e) ≈ x, (λ(e,\longunderscore).(\target␣(\widetilde{g}␣e),\longunderscore)))\),
for which we can check that \(π₀\) is a \(\T^∞\)-coalgebra homomorphism:
\begin{align*}
    \T^∞ π₀␣&(\tar_g␣(x,\longunderscore)) \\
        &= \left(∑_{e : \overline{g}} \source␣(\widetilde{g}␣e) ≈ x, π₀ ∘ (λ(e,\longunderscore).(\target␣(\widetilde{g}␣e),\longunderscore))\right) \\
        &= \left(∑_{e : \overline{g}} \source␣(\widetilde{g}␣e) ≈ x, (λ(e,\longunderscore).\target␣(\widetilde{g}␣e))\right) \\
        &= \left(∑_{y:V}〈x,y〉∈g ,π₀\right) \\
        &= m_g␣x \\
        &= m_g␣(π₀␣(x,\longunderscore)) \tag*{\qedhere}
\end{align*}\end{proof}

\textbf{Remark:} For \(\tar_g\), just as for \(m_g\), we will slightly
abuse notation, justified by Proposition \ref{prop:graph-coalg-Target},
and write:
\[\tar_g␣(x,\longunderscore) = \left(∑_{y:V}〈x,y〉∈g,λ(y,e).(y,|(x,e)|)\right).\]
Again, a more careful treatment is found in the formalisation.

\begin{lemma}[\agdalink{https://elisabeth.stenholm.one/non-wellfounded-set-theory/v4/e-structure.graph.to-P-n-coalgebra.html\#7669}]

\label{lem:truncation-of-n}If a graph \(g : V\) is an \(n\)-type in
\((V,∈)\) then \(\tar_g\) factors through a \(\Tⁿ\)-coalgebra
\(\tar_{n,g} : \Target g → \Tⁿ\Target g\).\end{lemma}

\begin{proof}

We only need to prove that \(π₁(\tar_g␣(x,\longunderscore))\) is
\((n-1)\) truncated, but with the notation we just introduced:
\(π₁(\tar_g␣(x,\longunderscore)):∑_{y:V}〈x,y〉∈g → \Target g\) is
defined by \(π₁(\tar_g␣(x,\longunderscore)) = λ(y,e).(y,|(x,e)|)\).
Since the map is the identity on the base, it suffices to consider the
truncation level of the second component:
\(〈x,y〉∈g → ∃_{x:V}〈x,y〉∈g\). The codomain is a proposition, so the
map's truncation level is the same as that of the domain, which by
assumption is \((n-1)\).\end{proof}

\subsection{Aczel's anti-foundation
axiom}\label{aczels-anti-foundation-axiom}

Aczel's anti-foundation axiom can now be generalised to any truncation
level. We will demonstrate that if one could construct terminal
coalgebras for the \(\Tⁿ\) functors, the resulting ∈-structures would
satisfy the generalised axiom.

\begin{definition}[\agdalink{https://elisabeth.stenholm.one/non-wellfounded-set-theory/v4/e-structure.property.aczel-anti-foundation.html\#563}]

An ∈-structure \((V,∈)\), with an ordered pairing structure, has
\textbf{Aczel \(n\)-anti-foundation} (\(n\)-AFA), for \(n : ℕ^∞_{0}\),
if for every graph \(g : V\) the type of \(n\)-decorations of \(g\)~is
contractible. Equivalently, this can be split into two parts:

\begin{itemize}
\tightlist
\item
  \(n\)-AFA₁: For every graph \(g:V\) the type of \(n\)-decorations of
  \(g\) is inhabited
\item
  \(n\)-AFA₂: For every graph \(g:V\) the type of \(n\)-decorations of
  \(g\) is a proposition.
\end{itemize}

\end{definition}

The classical AFA axiom is equivalent to Aczel 0-anti-foundation, since
0-decorations are the usual decorations, and contractible is the HoTT
way of saying ``exists unique''.

As decorations are \(\T^∞_n\)-coalgebra homomorphisms, and in particular
\(\Tⁿ\)-coalgebra homomorphisms in \(n\)-level ∈-structures, one type
that would model AFA is the terminal \(\Tⁿ\)-coalgebra.

\begin{theorem}[\agdalink{https://elisabeth.stenholm.one/non-wellfounded-set-theory/v4/e-structure.property.aczel-anti-foundation.from-terminal-coalgebra.html\#5900}]

\label{thm:AFA-from-terminal-coalg}Suppose \((V,m)\) is the terminal
\(\Tⁿ\)-coalgebra and that \(V\) is locally \(U\)-small. Then the
induced ∈-structure has Aczel \(n\)-anti-foundation.\end{theorem}

\begin{proof}

It was shown in Gylterud and Stenholm~\cite{GylterudStenholm2026} (Theorem 1 and Theorem 16)
that \((V,m)\) has an ordered pairing structure. Let \(g : V\) be a
graph. By Proposition \ref{prop:equiv-decoration-Tn-coalgebra} we need
to show that the type of \(\T^∞_n\)-coalgebra homomorphisms from the
corresponding graph coalgebra \(m_g\), given by Proposition
\ref{prop:graph-coalg}, into \((V,m)\) is contractible. For these
propositions we need \((V,m)\) to be \(U\)-like, but this follows from
Lemma \ref{coalgebra-structures}. By Lemma \ref{lma:mg-factor} and the
fact that the map \(\Tⁿ␣V ↪ \T^∞_n␣V\) is an embedding and thus a
monomorphism, it is enough to show that the type of \(\Tⁿ\)-coalgebra
homomorphisms from \(m_{n,g}\) to \((V,m)\) is contractible. But this
follows from terminality of \((V,m)\).\end{proof}

\textbf{Remark:} In the proof above, if \((V,∈_m)\) has level \(n\), we
only use terminality with respect to coalgebras on \(n\)-types. Hence,
for an ∈-structure of level 0 it is sufficient to show terminality with
respect to mere sets.

\subsection{\texorpdfstring{Scott's anti-foundation axiom
\label{subsec:safa}}{Scott's anti-foundation axiom }}\label{scotts-anti-foundation-axiom}

Recall that, classically, SAFA is the statement that every Scott
extensional graph has a unique injective decoration and \(V\) itself is
Scott extensional. A graph is defined as being Scott extensional if
equality on the nodes is given by a tree isomorphism of the
corresponding unfolding trees. Note that two trees are isomorphic if
there is an isomorphism between the children of the roots, such that the
subtrees of each related pair of children are tree isomorphic. We can
see this as the unfolding step in a \(\T^∞\)-bismulation.

The terminal \(\T^∞\)-coalgebra, \(V^∞_∞\), can be thought of as the
type of trees, and the map induced by its terminality,
\(\corec^∞␣(A,m) : A → V^∞_∞\), is the unfolding of a coalgebra or graph
into a tree (starting in a given node). Because of univalence, the
identity type in \(V^∞_∞\) is equivalent to tree isomorphism. This means
that we can express Scott extensionality for a graph as saying that
\(\corec^∞␣(\Target␣g,\tar_g)\) is an embedding. Every function in HoTT
has an associated action on paths, which becomes an equivalence for an
embedding. So, if \(\corec^∞␣(\Target␣g,\tar_g)\) is an embedding, its
action on paths of the graph provides an equivalence between equality in
the graph and isomorphism of its unfolding trees.

On higher type levels, it is a bit strong to require an embedding. For
instance, in multisets (which are the material set theory equivalent of
groupoids), we would like to consider a graph like
\(\{〈∅,∅〉,〈∅,∅〉\}₁\) as a Scott extensional representation of the
complete binary tree. However, this tree has many non-trivial
automorphisms in \(V^∞_∞\), which our single node, ∅, does not have. An
embedding would require nodes in the graph to come prefilled with these
automorphisms, but in our models this is not required. However, some
restriction must be enforced to make the type levels of the left and
right hand side of the equivalence defining an ∞-decoration agree: see
the equivalence (\ref{eq:infinity-decoration}) just after Definition
\ref{def:n-decoration}. We therefore define the notion of a graph being
Scott \(n\)-extensional as follows.

\begin{definition}[\agdalink{https://elisabeth.stenholm.one/non-wellfounded-set-theory/v4/e-structure.property.scott-anti-foundation.html\#1690}]

Given a graph \(g:V\) and \(n : ℕ^∞_{-1}\), we say that \(g\) is Scott
\((n+1)\)-extensional if the tree unfolding map
\(\corec^∞␣(\Target␣g,\tar_g)\) is \(n\)-truncated.\end{definition}

Clearly, being Scott \(n\)-extensional implies being Scott
\((n+1)\)-extensional, and by the reasoning above, Scott
\(0\)-extensionality is the usual notion of Scott extensionality in
level 0 ∈-structures. Furthermore, if the graph is a set level graph
(meaning that \(\Target␣g\) is a set and \(\tar_g\) factors through
\(\T¹\)), then it is automatically Scott \(1\)-extensional.

We can now define Scott's anti-foundation axiom for ∈-structures of any
type level.

\begin{definition}[\agdalink{https://elisabeth.stenholm.one/non-wellfounded-set-theory/v4/e-structure.property.scott-anti-foundation.html\#2232}]

A \(U\)-like ∈-structure \((V,∈)\), with an ordered pairing structure,
satisfies \textbf{Scott \(n\)-anti-foundation} (\(n\)-SAFA), for
\(n : ℕ^∞₀\), if the two properties \(n\)-SAFA₁ and SAFA₂ hold:

\begin{itemize}
\tightlist
\item
  \(n\)-SAFA₁: Any Scott \(n\)-extensional graph \(g:V\) has an
  ∞-decoration.
\item
  SAFA₂: For any graph \(g\) the type of ∞-decorations is a proposition.
\end{itemize}

\end{definition}

The classical notion of SAFA then corresponds to what is defined above
as Scott 0-anti-foundation. SAFA₂ is the same as ∞-AFA₂, and since being
Scott ∞-extensional is a vacuous requirement, we get that ∞-SAFA is
equivalent to ∞-AFA.

\section{The coiterative hierarchy}\label{the-coiterative-hierarchy}

The coiterative hierarchy is a dualisation of a specific construction of
the iterative hierarchy \citep{gylterud-iterative}. That construction
starts with the type of all \emph{wellfounded} trees and picks out the
subset of those which are hereditarily sets (i.e.~in each node each
immediate subtree is unique). Figure \ref{fig:tree1} and \ref{fig:tree2}
give examples of which trees are and which are not heriditary sets. The
coiterative hierarchy is constructed dually, starting from the type of
all (possibly non-wellfounded) trees, and picking out those which are
co-hereditarily sets. That is, no matter how far we go into the tree, in
each node the immediate subtrees are always distinct. Figure
\ref{fig:infinite-binary-tree} and \ref{fig:infinite-binary-tree-2}
gives examples to illustrate the notion of trees being co-hereditary
sets.

\begin{figure}[ht!]
    \centering
    \begin{minipage}{0.45\textwidth}
        \centering
        \begin{tikzpicture}[level distance=1.5cm, sibling distance=3cm, scale=0.6]
            \node {•} 
                  child {
                    node {•}
                    child { node {•}
                            child { node {•} } }
                    child { node {•} }
                  }
                  child { node {•} }
                ;
        \end{tikzpicture}
        \caption{This tree represents an iterative set, namely: $\{\{\{∅\}₀,∅\}₀,∅\}₀$.}
        \label{fig:tree1}
    \end{minipage}\hfill
    \begin{minipage}{0.45\textwidth}
        \centering
        \begin{tikzpicture}[level distance=1.5cm, sibling distance=3cm, scale=0.6]
            \node {•}
                  child {
                    node {•}
                    child { node {•} }
                    child { node {•} }
                  }
                  child { node {•} }
                ;
        \end{tikzpicture}
        \caption{This tree does not represent an iterative set
                 because the left child of the root has two equal
                 children. It does however represent the iterative multiset: $\{\{∅,∅\}₁,∅\}₁$.}
        \label{fig:tree2}
    \end{minipage}
\end{figure}

\begin{figure}[ht!]
    \centering
    \begin{minipage}{0.45\textwidth}
    \centering
    \begin{tikzpicture}[level distance=1.5cm,level/.style={sibling distance={5cm/max(1,#1)}}, grow=south, scale=0.6]
        \node {•}
            child { node {•}
                child { node {•}
                    child { node {} edge from parent[dotted]
                        child[missing]
                        child[missing]
                    }
                    child { node {} edge from parent[dotted]
                        child[missing]
                        child[missing]
                    }
                }
                child { node {•}
                    child { node {} edge from parent[dotted]
                        child[missing]
                        child[missing]
                    }
                    child { node {} edge from parent[dotted]
                        child[missing]
                        child[missing]
                    }
                }
            }
            child { node {•}
                child { node {•}
                    child { node {} edge from parent[dotted]
                        child[missing]
                        child[missing]
                    }
                    child { node {} edge from parent[dotted]
                        child[missing]
                        child[missing]
                    }
                }
                child { node {•}
                    child { node {} edge from parent[dotted]
                        child[missing]
                        child[missing]
                    }
                    child { node {} edge from parent[dotted]
                        child[missing]
                        child[missing]
                    }
                }
            };
        \node at (1.2, -4.5) {\vdots};
        \node at (3.75, -4.5) {\vdots};
        \node at (-1.2, -4.5) {\vdots};
        \node at (-3.75, -4.5) {\vdots};
    \end{tikzpicture}
    \caption{The full binary tree is not a coiterative set. But rather a multiset $b =\{b,b\}₁$.}
    \label{fig:infinite-binary-tree}
\end{minipage}\hfill
\begin{minipage}{0.45\textwidth}
    \centering
    \begin{tikzpicture}[level distance=1.5cm,level/.style={sibling distance={6cm/max(1.5,#1)}}, grow=south, scale=0.6]
        \node{x}
        child{ node{x}
                child{  node{\textcolor{gray}{x}}
                        child {node {} edge from parent[dotted]}
                        child {node {} edge from parent[dotted]}}
                child{  node{\textcolor{gray}{y}}
                        child {node {} edge from parent[dotted]}}}
        child{ node{y}
                child{ node{x}
                        child{  node{\textcolor{gray}{x}}
                                child {node {} edge from parent[dotted]}
                                child {node {} edge from parent[dotted]}}
                        child{  node{\textcolor{gray}{y}}
                                child {node {} edge from parent[dotted]}}
                    }}
        ;
    \end{tikzpicture}
    \caption{This infinite binary tree represents the coiterative set $x$ which is part of the solution to the equations $x=\{x,y\}₀$ and $y=\{x\}₀$}
    \label{fig:infinite-binary-tree-2}
\end{minipage}
\end{figure}

\subsection{\texorpdfstring{Defining
\(Vⁿ_∞\)}{Defining \textbackslash Vⁿ\_∞}}\label{defining-vux207f_}

In \emph{Univalent Material Set Theory} \citep{GylterudStenholm2026},
the construction of an iterative hierarchy of sets was extended to a
hierarchy of \(n\)-types, \(Vⁿ\). When dualising to coiterative sets we
will keep this level of generality and construct a coiterative hierarchy
of \(n\)-types, \(Vⁿ_∞\). The first level, \(V⁰_∞\), yields the type of
coiterative sets.

The iterative hierarchy was carved out from the W-type
\(\V^∞ ≔ \W{A : U}{A}\), as a subtype, using an inductive predicate
\(\ittype{n} : \V^∞ → \Type\). The coiterative hierarchy will, dually,
be carved out as a subtype from the M-type, \(\V^∞_∞ ≔ \M{A : U}{A}\),
and a coinductive predicate \(\iscoittype{n} : \V^∞_∞ → \Type\). But
first we need an auxiliary predicate on each depth.

\begin{definition}[\agdalink{https://elisabeth.stenholm.one/non-wellfounded-set-theory/v4/coiterative.set.html\#2088}]

For \(n:\Nat_{-1}\), define the predicate: \begin{align*}
    &\iscoittype{(n+1)} : ℕ → \V^∞_∞ → \Type \\
    &\iscoittypen{(n+1)}{0}␣x ≔ \isntruncmap{n}␣\widetilde{x}\\
    &\iscoittypen{(n+1)}{(k+1)}␣x ≔ ∏_{a : \overline{x}} \iscoittypen{(n+1)}{k} \left(\widetilde{x}␣a\right).
\end{align*}\end{definition}

\begin{definition}[\agdalink{https://elisabeth.stenholm.one/non-wellfounded-set-theory/v4/coiterative.set.html\#2934}]

For \(n:\Nat\), define the predicate: \begin{align*}
    &\iscoittype{n} : \V^∞_∞ → \Type \\
    &\iscoittype{n}␣x ≔ ∏_{k : ℕ} \iscoittypen{n}{k}␣x
\end{align*}\end{definition}

The predicate \(\iscoittype{n}\) is clearly a proposition, and we now
define the type of coiterative \(n\)-types, as the resulting subtype of
\(\V^∞_∞\):

\begin{definition}[The coiterative hierarchy \agdalink{https://elisabeth.stenholm.one/non-wellfounded-set-theory/v4/coiterative.set.html\#3554}]

For \(n : \Nat\), let \(\Vⁿ_∞\) denote the type of coiterative
\(n\)-types: \[\Vⁿ_∞ ≔ ∑_{x:\V^∞_∞}\iscoittype{n}␣x.\]\end{definition}

\begin{proposition}[\agdalink{https://elisabeth.stenholm.one/non-wellfounded-set-theory/v4/coiterative.set.html\#4717}]

\label{prop:Vn-is-subtype-of-Vinf}\(\Vⁿ_∞\) is a subtype of \(\V^∞_∞\),
i.e. the projection \(\Vⁿ_∞ → \V^∞_∞\) is an embedding.\end{proposition}

This means, in particular, that the identity type on \(\Vⁿ_∞\) is the
same as the identity type on \(\V^∞_∞\).

One of the requirements to satisfy SAFA is that the type of
∞-decorations is a proposition. By the characterisation of ∞-decorations
as \(\T^∞\)-coalgebra homomorphisms it is sufficient for the model to be
a simple \(\T^∞\)-coalgebra. Thus, we show this for \(\Vⁿ_∞\).

\begin{proposition}[\agdalink{https://elisabeth.stenholm.one/non-wellfounded-set-theory/v4/coiterative.set.html\#9418}]

\label{prop:vnif-simple}\((\Vⁿ_∞,\desupⁿ)\) is a simple
\(\T^∞\)-coalgebra.\end{proposition}

\begin{proof}

Let \((X,m)\) be a \(\T^∞\)-coalgebra. Since \(\Vⁿ_∞\) embeds into
\(\V^∞_∞\) by Proposition \ref{prop:Vn-is-subtype-of-Vinf} and since
polynomial functors preserve embeddings, it follows that the
type\newline \(\Hom{\Coalg{\T^∞}}␣(X,m)␣(\Vⁿ_∞, \desupⁿ)\) embeds into
the type \(\Hom{\Coalg{\T^∞}}␣(X,m)␣(\V^∞_∞, \desup^∞)\). The latter
type is contractible and thus a proposition. The result then follows
from the fact that any type which embeds into a proposition is a
proposition.\end{proof}

\subsection{\texorpdfstring{\(\Vⁿ_∞\) is a fixed point for
\(\Tⁿ\)}{\textbackslash Vⁿ\_∞ is a fixed point for \textbackslash Tⁿ}}\label{vux207f_-is-a-fixed-point-for-tux207f}

The elements in \(\Vⁿ_∞\) are non-wellfounded trees where all branchings
are \((n-1)\)-truncated maps. So when one removes the root from a tree,
one gets a small type and an \((n-1)\)-truncated map from that type into
\(\Vⁿ_∞\). Similarly, if one has a small type and an \((n-1)\)-truncated
map from that type into \(\Vⁿ_∞\) then one can construct a tree in
\(\Vⁿ_∞\) by adding a root node. Hence, we will show that \(\Vⁿ_∞\) is a
fixed point of \(\Tⁿ\).

\begin{lemma}[\agdalink{https://elisabeth.stenholm.one/non-wellfounded-set-theory/v4/coiterative.set.html\#3755}]

\label{lma:equiv-iscoittype}For any \(x : \V^∞_∞\), there is an
equivalence \[\iscoittype{n}␣x ≃ \left(\isntruncmap{n}␣\widetilde{x} × 
    ∏_{a : \overline{x}} \iscoittype{n}
    \left(\widetilde{x}␣a\right)\right).\]\end{lemma}

\begin{proof}

Follows by induction over ℕ.\end{proof}

\begin{theorem}[\agdalink{https://elisabeth.stenholm.one/non-wellfounded-set-theory/v4/coiterative.set.html\#7701}]

\label{thm:Vn-fixed-point}\(\Vⁿ_∞\) is a fixed point for
\(\Tⁿ\).\end{theorem}

\begin{proof}

We begin with the case \(n = ∞\). Since \(\V^∞_∞\) is the terminal
\(\T^∞\)-coalgebra, it is in particular a fixed point for \(\T^∞\).

For the case \(n < ∞\), let \(x : \Vⁿ_∞\), then by Lemma
\ref{lma:equiv-iscoittype}, the element \((\overline{x},\widetilde{x})\)
lies in \(\Tⁿ␣\Vⁿ_∞\). By the same token, given \(A : U\) and
\(f : A ↪_{n-1} \Vⁿ_∞\), the element \(\ssup^∞␣(A , f)\) is a
coiterative \(n\)-type.\end{proof}

For the two maps given by Theorem \ref{thm:Vn-fixed-point} we introduce
the following notation: \begin{align*}
    &\desupⁿ : \Vⁿ_∞ → \Tⁿ␣\Vⁿ_∞,\\
    &\operatorname{sup}ₙ : \Tⁿ␣\Vⁿ_∞ → \Vⁿ_∞.
\end{align*}

\begin{proposition}[\agdalink{https://elisabeth.stenholm.one/non-wellfounded-set-theory/v4/coiterative.set.html\#8437}]

\label{prop:inclusion-Vn-is-hom}The inclusion \(\Vⁿ_∞ ↪ \V^∞_∞\) is a
\(\T^∞\)-coalgebra homomorphism from \((\Vⁿ_∞, \desupⁿ)\) (seen as a
\(\T^∞\)-coalgebra) to \((\V^∞_∞, \desup^∞)\).\end{proposition}

\begin{proof}

This holds definitionally, i.e.~the homotopy is given by
\(\reflhtpy\).\end{proof}

\subsection{\texorpdfstring{Non-terminality of \(\V⁰_∞\) as a
\(\T⁰\)-coalgebra}{Non-terminality of \textbackslash V⁰\_∞ as a \textbackslash T⁰-coalgebra}}\label{non-terminality-of-vux2070_-as-a-tux2070-coalgebra}

Even though \(\Vⁿ_∞\) is a fixed point for \(\Tⁿ\) and is a subtype of
the terminal \(\T^∞\)-coalgebra, it turns out \textbf{not} to be the
terminal \(\Tⁿ\)-coalgebra. At least \(\V⁰_∞\) is not the terminal
\(\T⁰\)-coalgebra. But we conjecture this result to hold for all \(n\).
This is surprising since the dual construction gives the initial algebra
of \(\Tⁿ\) (Theorem 15 of Gylterud and Stenholm~\cite{GylterudStenholm2026}). Intuitively,
the reason is that in the wellfounded setting tree isomorphism coincides
with bisimulation, while in the non-wellfounded setting it does not. In
Subsection \ref{the-coiterative-hierarchy-as-a-model-of-set-theory} we
shall see how this means that \(\Vⁿ_∞\) satisfy (the generalisations of)
Scott's anti-foundation axiom rather than Aczel's.

For \(\Vⁿ_∞\) to be terminal, any graph (considered as a
\(\Tⁿ\)-coalgebra) should have a unique representative in \(\Vⁿ_∞\). But
\(\Vⁿ_∞\) contains more than one representative of some graphs, i.e.~we
can construct a \(\Tⁿ\)-coalgebra for which there are two distinct
\(\Tⁿ\)-coalgebra homomorphisms into \(\Vⁿ_∞\). One of the maps sends
each node to its unfolding tree. Because the functorial action of
\(\Tⁿ\) takes the \((n-1)\)-image of the composite map, i.e.~it
collapses some structure, there is also a \(\Tⁿ\)-coalgebra homomorphism
which maps the nodes to different trees.

\begin{theorem}

\label{thm:non-terminality}\(\V⁰_∞\) is not the terminal
\(\T⁰\)-coalgebra.\end{theorem}

\begin{proof}

Consider the following \(\T⁰\)-coalgebra \((X,m)\), represented as a
graph:

\begin{center}
    \begin{tikzcd}
        x \arrow[r, bend left] \arrow[loop, distance=2em, in=215, out=145] & y \arrow[l, bend left]
    \end{tikzcd}
\end{center}

The unfolding trees of the two nodes as given by
\(\corec^∞ : X → \V^∞_∞\) are distinct (\(\corec^∞␣x\) being the tree
depicted in Figure \ref{fig:infinite-binary-tree-2}), so \(\corec^∞\)
factors as a \(\T⁰\)-coalgebra homomorphism, \(f : X → \V⁰_∞\), from
\((X,m)\) to \((\V⁰_∞,
\desup₀)\), such that \(f␣x ≠ f␣y\).

On the other hand, let \(g\) be the map that sends both nodes to the
infinite unary tree, which we will denote \(q : V⁰_∞\):

\begin{center}
    \begin{tikzcd}
        \bullet \arrow[d] \\
        \bullet \arrow[d] \\
        \vdots           
    \end{tikzcd}
\end{center}

Clearly, \(g\) is also a \(\T⁰\)-coalgebra homomorphism:
\[\T⁰g␣(m␣x) = (\image␣(g ∘ \widetilde{x}), \incl) = (1,λ\longunderscore.q) = (\overline{q},\widetilde{q}) = \desup₀(g␣x)\]
and likewise for \(y\):
\[\T⁰g␣(m␣y) = (\image␣(g ∘ \widetilde{y}), \incl) = (1,λ\longunderscore.q) = (\overline{q},\widetilde{q}) = \desup₀(g␣y).\]

However, since \(f␣x≠f␣y\) and \(g␣x=g␣y\), we get that \(f\) and \(g\)
are two distinct \(\T⁰\)-coalgebra homomorphisms from \((X,m)\) to
\((\V⁰_∞,\desup₀)\).\end{proof}

Note that the result above does not contradict the fact that \(\Vⁿ_∞\)
is simple as a \(\T^∞\)-coalgebra. The proof of non-terminality of
\(\V⁰_∞\) demonstrates that it is not simple as a \(\T⁰\)-coalgebra.

\subsection{\texorpdfstring{Local smallness of
\(\Vⁿ_∞\)}{Local smallness of \textbackslash Vⁿ\_∞}}\label{local-smallness-of-vux207f_}

The functorial action of \(\Tⁿ\) takes the \((n-1)\)-image of a map. In
order for this to be small, the domain must be small and the codomain
appropriately locally small. In particular, when we are considering maps
into \(\Vⁿ_∞\), we use the fact that this type is locally small, as we
will show in this section. This result uses univalence and follows from
the characterisation of the identity on an M-type as an indexed M-type.

The idea is that, by univalence, the indexed functor
\(\E_{(V^∞_∞,\desup^∞)}\) (Definition \ref{def:E-functor}) is equivalent
to the indexed functor \(\E^≃_{(V^∞_∞,\desup^∞)}\), for which the
corresponding indexed M-type is small.

\begin{definition}[\agdalink{https://elisabeth.stenholm.one/non-wellfounded-set-theory/v4/functor.slice.html\#3159}]

Given \(X : \Type\) and \(m : X → \left(∑_{A : \Type} A → X\right)\), we
define the \(\left(X × X\right)\)-indexed polynomial functor
\begin{align*}
    &\E^≃_{(X , m)} : \left(X × X → \Type\right) → \left(X × X → \Type\right)\\
    &\E^≃_{(X , m)}␣R\;(x,y) ≔ ∑_{e : \overline{x}≃\overline{y}} ∏_{a : \overline{x}}
        R␣(\widetilde{x}␣a, \widetilde{y}␣(e␣a)).
 \end{align*} The functorial action is postcomposition on the second
component. The functor respects identities and composition
definitionally.\end{definition}

The subscript \((X,m)\) might be omitted if it is clear from the
context.

\begin{definition}

Given \(X : \Type\) and \(m : X → \left(∑_{A : \Type} A → X\right)\), we
define the following types:

\begin{itemize}
\tightlist
\item
  The type of \(\E^≃_{(X,m)}\)-coalgebras is
  \[\Coalg{\E^≃_{(X,m)}} ≔ ∑_{R : X × X → \Type} ∏_{(x,y) : X × X} R␣(x,y) → \E^≃␣R␣(x,y).\]
\item
  Given two \(\E^≃_{(X,m)}\)-coalgebras \((R,α)\) and \((Q,β)\), the
  type of \(\E^≃_{(X,m)}\)-coalgebra homomorphisms between them is
  \begin{align*}
    &\Hom{\Coalg{\E^≃_{(X,m)}}}␣(R,α)␣(Q,β) \\
    &\hspace{20pt}≔ ∑_{f : ∏_{(x,y) : X × X} R␣(x,y) → Q␣(x,y)} 
        ∏_{(x,y) : X × X} β␣(x,y) ∘ f␣(x,y) ∼ \E^≃␣f␣(x,y) ∘ α␣(x,y).
  \end{align*}
\end{itemize}

\end{definition}

The identity type is an \(\E^≃_{(X,m)}\)-coalgebra, for any pair
\((X,m)\).

\begin{definition}[\agdalink{https://elisabeth.stenholm.one/non-wellfounded-set-theory/v4/functor.slice.html\#4042}]

Define the following map by path induction: \begin{align*}
    &\idecoalg_{(X,m)}^≃ : ∏_{(x,y) : X × X} x = y → \E^≃_{(X,m)}␣({=})\;(x,y) \\
    &\idecoalg_{(X,m)}^≃␣(x,x)␣\refl ≔ (\idequiv, \reflhtpy).
\end{align*} The pair \((=,\idecoalg_{(X,m)}^≃)\) is the
\textbf{identity \(\E^≃_{(X,m)}\)-coalgebra}.\end{definition}

There is an equivalence between \(\E\)-coalgebras and
\(\E^≃\)-coalgebras.

\begin{proposition}[\agdalink{https://elisabeth.stenholm.one/non-wellfounded-set-theory/v4/coiterative.infinity-multiset.html\#2307}]

\label{prop:E-E'-naturally-equiv}Given \(X : \Type\) and
\(m : X → \left(∑_{A :
\Type} A → X\right)\), there is a natural family of equivalences
\[\E_{(X , m)}␣R\;(x,y) ≃ \E^≃_{(X , m)}␣R\;(x,y).\]\end{proposition}

\begin{proof}

Follows by univalence.\end{proof}

This gives us an alternative characterisation of the identity type on
\(\V^∞_∞\).

\begin{theorem}[\agdalink{https://elisabeth.stenholm.one/non-wellfounded-set-theory/v4/coiterative.infinity-multiset.html\#5783}]

\label{thm:identity-on-Vinf}The identity coalgebra \((=,\idecoalg^≃)\)
is the terminal \(\E^≃_{(\V^∞_∞ , \desup^∞)}\)-coalgebra.\end{theorem}

\begin{proof}

By Theorem \ref{thm:identity-on-m-types}, the identity coalgebra
\((=,\idecoalg)\) is the terminal
\(\E_{(\V^∞_∞ , \desup^∞)}\)-coalgebra. Since the functors \(\E\) and
\(\E^≃\) are naturally equivalent by Proposition
\ref{prop:E-E'-naturally-equiv}, and this equivalence maps
\((=,\idecoalg)\) to \((=,\idecoalg^≃)\), the identity coalgebra
\((=,\idecoalg^≃)\) is the terminal coalgebra for
\(\E^≃_{(\V^∞_∞ , \desup^∞)}\).\end{proof}

Note that by the theorem above, there is for any \(x,y : \V^∞_∞\) an
equivalence
\[(x = y) ≃ ∑_{e : \overline{x}≃\overline{y}} ∏_{a : \overline{x}}
        \widetilde{x}␣a = \widetilde{y}␣(e␣a).\]

\begin{theorem}[\agdalink{https://elisabeth.stenholm.one/non-wellfounded-set-theory/v4/coiterative.infinity-multiset.html\#6120}]

\label{thm:Vinf-is-locally-small}\(\V^∞_∞\) is locally
\(U\)-small.\end{theorem}

\begin{proof}

Since \(\E^≃_{(\V^∞_∞ , \desup^∞)}\) is an indexed polynomial functor,
it has a corresponding indexed M-type which is the terminal
\(\E^≃_{(\V^∞_∞ ,
\desup^∞)}\)-coalgebra. In their paper on non-wellfounded trees in HoTT,
Ahrens et al.~\cite{ahrens-2015} construct indexed M-types from inductive types. They
only spell out the details of the non-indexed case, leaving the details
of the indexed case in the accompanying formalisation. They note that
the construction of the non-indexed M-type does not raise the universe
level, in the proof of Theorem 7. Though they do not state so
explicitly, this is also the case for the indexed M-types. In
particular, the carrier of the indexed M-type is the limit of successive
applications of the indexed polynomial functor in question, to the unit
type. As the universe level of an indexed polynomial functor applied to
a type and an index does not depend on the level of the indexing type,
the corresponding indexed M-type also does not depend on the level of
the indexing type. (One can also see that this is the case in the
formalisation.) In our case, the universe level of the indexed M-type
corresponding to \(\E^≃_{(\V^∞_∞ , \desup^∞)}\) is the least upper bound
of the universe levels of \(\overline{x}≃\overline{y}\) and
\(\overline{x}\), which is \(U\).

Since (the carriers of) any two terminal \(\E^≃_{(\V^∞_∞ ,
\desup^∞)}\)-coalgebras are equivalent, it follows that \(\V^∞_∞\) is
locally \(U\)-small.\end{proof}

\begin{corollary}[\agdalink{https://elisabeth.stenholm.one/non-wellfounded-set-theory/v4/coiterative.set.html\#7949}]

\label{cor:Vn-is-locally-small}\(\Vⁿ_∞\) is locally
\(U\)-small.\end{corollary}

\begin{proof}

By Proposition \ref{prop:Vn-is-subtype-of-Vinf}, \(\Vⁿ_∞\) is a subtype
of \(\V^∞_∞\) and thus has the same identity type. The result then
follows from the fact that \(\V^∞_∞\) is locally \(U\)-small, by Theorem
\ref{thm:Vinf-is-locally-small}.\end{proof}

\subsection{\texorpdfstring{Coalgebra homomorphisms into
\(\Vⁿ_∞\)}{Coalgebra homomorphisms into \textbackslash Vⁿ\_∞}}\label{coalgebra-homomorphisms-into-vux207f_}

How do we construct a map from a \(\Tⁿ\)-coalgebra, say \((X , m)\),
into \(\Vⁿ_∞\)? An obvious approach is to view \((X , m)\) as a
\(\T^∞\)-coalgebra and show that \(\corec^∞ : X → \V^∞_∞\) lands in
\(\Vⁿ_∞\), where \(\corec^∞\) is the underlying map of the unique
\(\T^∞\)-coalgebra homomorphism from \((X,m)\) to \((\V^∞_∞,\desup^∞)\).
Unfortunately, \(\corec^∞\) does not always land in \(\Vⁿ_∞\).

Viewing \((X,m)\) as a graph, \(\corec^∞\) maps each node to its
unfolding tree. Consider now the \(\T⁰\)-coalgebra represented by the
following graph:

\begin{center}
    \begin{tikzcd}
                                                        & \bullet \arrow[rd] \arrow[ld] &                                                     \\
    \bullet \arrow[loop, distance=2em, in=305, out=235] &                               & \bullet \arrow[loop, distance=2em, in=305, out=235]
    \end{tikzcd}
\end{center}

The topmost node is mapped by \(\corec^∞\) to the tree

\begin{center}
    \begin{tikzcd}
                      & \bullet \arrow[rd] \arrow[ld] &                   \\
    \bullet \arrow[d] &                               & \bullet \arrow[d] \\
    \bullet \arrow[d] &                               & \bullet \arrow[d] \\
    \vdots            &                               & \vdots           
    \end{tikzcd}
\end{center}

which is not an element of \(\V⁰_∞\) as the branching at the root is not
an embedding.

However, if \(\corec^∞\) is an \((n-1)\)-truncated map, then it lands in
\(\Vⁿ_∞\).

\begin{proposition}[\agdalink{https://elisabeth.stenholm.one/non-wellfounded-set-theory/v4/coiterative.set.html\#18240}]

\label{prop:corec-lands-in-Vn}Given a \(\Tⁿ\)-coalgebra \((X,m)\), if
\(\corec^∞ : X → \V^∞_∞\) is an \((n-1)\)-truncated map, then for all
\(x : X\), \(\corec^∞␣x\) is a coiterative \(n\)-type.\end{proposition}

\begin{proof}

For \(x : X\) we need to show that
\[∏_{k:ℕ} \iscoittypen{n}{k}␣(\corec^∞␣x).\] Proceed by induction on
\(k\).

For the base case, note that since \(\corec^∞\) is a \(\T^∞\)-coalgebra
homomorphism, we have
\[\widetilde{(\corec^∞␣x)} = \corec^∞ ∘\,\widetilde{x}.\] Both these
maps are \((n-1)\)-truncated, and therefore the composition is
\((n-1)\)-truncated.

Similarly, for the induction step, since \(\corec^∞\) is a homomorphism,
it is enough to show that
\[∏_{a : \overline{x}} \iscoittypen{n}{k}␣(\corec^∞␣(\widetilde{x}␣a)).\]
But this follows from the induction hypothesis.\end{proof}

\begin{definition}[\agdalink{https://elisabeth.stenholm.one/non-wellfounded-set-theory/v4/coiterative.set.html\#18833}]

\label{corec-n}Given a \(\Tⁿ\)-coalgebra \((X,m)\) for which
\(\corec^∞\) is an \((n-1)\)-truncated map, let \[\corecⁿ : X → \Vⁿ_∞\]
denote the restriction of \(\corec^∞\) into \(\Vⁿ_∞\) by Proposition
\ref{prop:corec-lands-in-Vn}.\end{definition}

The map \(\corecⁿ\) is a \(\Tⁿ\)-coalgebra homomorphism. This is an
instance of a useful lemma about which maps into \(\Vⁿ_∞\) are
\(\Tⁿ\)-coalgebra homomorphisms.

\begin{lemma}[\agdalink{https://elisabeth.stenholm.one/non-wellfounded-set-theory/v4/coiterative.set.html\#14232}]

\label{lma:equiv-Tn-Tinf-homomorphism}Let \((X,m)\) be a
\(\Tⁿ\)-coalgebra and let \(f : X ↪_{n-1} \Vⁿ_∞\). Then there is an
equivalence of types between \(f\) being a \(\Tⁿ\)-coalgebra
homomorphism and \(π₀ ∘ f\) being a \(\T^∞\)-coalgebra
homomorphism.\end{lemma}

\begin{proof}

Consider the following diagram:

\begin{center}
    \begin{tikzcd}
        X \arrow[d, "m"'] \arrow[rr, "f"]       &  & \Vⁿ_∞ \arrow[d, "\desupⁿ"] \arrow[rr, hook, "π₀"] &  & \V^∞_∞ \arrow[dd, "\desup^∞"] \\
        \Tⁿ␣X \arrow[d, hook] \arrow[rr, "\Tⁿ␣f"] &  & \Tⁿ␣\Vⁿ_∞ \arrow[d, hook]                       &  &                               \\
        \T^∞␣X \arrow[rr, "\T^∞␣f"']        &  & \T^∞␣\Vⁿ_∞ \arrow[rr, hook, "\T^∞␣π₀"']       &  & \T^∞␣\V^∞_∞                
    \end{tikzcd}
\end{center}

The maps from \(\Tⁿ␣X\) to \(\T^∞␣X\) and from \(\Tⁿ␣\Vⁿ_∞\) to
\(\T^∞␣\Vⁿ_∞\) are embeddings as they simply forgets that the map in the
second coordinate is \((n-1)\)-truncated. Additionally, \(\T^∞␣π₀\) is
an embedding since π₀ is an embedding and postcomposition by an
embedding is an embedding.

The square on the right commutes as the inclusion of \(\Vⁿ_∞\) into
\(\V^∞_∞\) is a \(\T^∞\)-coalgebra homomorphism (Proposition
\ref{prop:inclusion-Vn-is-hom}), and since \(f\) is an
\((n-1)\)-truncated map, the lower left square also commutes. Note that
commutativity of the top left square means that \(f\) is a
\(\Tⁿ\)-coalgebra homomorphism, and that commutativity of the outermost
square means that \(π₀ ∘ f\) is a \(\T^∞\)-coalgebra homomorphism.

Since embeddings are monomorphisms, a filling of the upper left square
is equivalent to an equality between two maps in question when
postcomposed with the forgetful map and \(\T^∞␣π₀\). By the
commutativity of the right square and the lower left square, this is in
turn equivalent to a filling of the outer square.\end{proof}

\begin{proposition}[\agdalink{https://elisabeth.stenholm.one/non-wellfounded-set-theory/v4/coiterative.set.html\#18951}]

\label{prop:corec-n-is-Tn-homomorphism}Let \((X,m)\) be a
\(\Tⁿ\)-coalgebra for which \(\corec^∞\) is an \((n-1)\)-truncated map,
then \(\corecⁿ\) is an \((n-1)\)-truncated map, and it is a
\(\Tⁿ\)-coalgebra homomorphism into
\((\Vⁿ_∞,\desupⁿ)\).\end{proposition}

\begin{proof}

Recall that \[π₀ ∘ \corecⁿ ≡ \corec^∞.\] Since \(\corec^∞\) is an
\((n-1)\)-truncated map and \(π₀ : \Vⁿ_∞ → \V^∞_∞\) is an embedding, it
follows that \(\corecⁿ\) is an \((n-1)\)-truncated-map. By Lemma
\ref{lma:equiv-Tn-Tinf-homomorphism}, since \(\corec^∞\) is a
\(\T^∞\)-coalgebra homomorphism, it follows that \(\corecⁿ\) is a
\(\Tⁿ\)-coalgebra homomorphism.\end{proof}

Even though \((\Vⁿ_∞,\desupⁿ)\) is not the terminal \(\Tⁿ\)-coalgebra,
it is \emph{almost} terminal---it is terminal with respect to truncated
maps.

\begin{theorem}[\agdalink{https://elisabeth.stenholm.one/non-wellfounded-set-theory/v4/coiterative.set.html\#20255}]

\label{thm:Vn-terminal-wrt-truncated-maps}Let \((X,m)\) be a
\(\Tⁿ\)-coalgebra for which \(\corec^∞\) is an \((n-1)\)-truncated map.
Then the following type is contractible:
\[∑_{(f,α) : \Hom{\Coalg{\Tⁿ}}␣(X,m)␣(\Vⁿ_∞,\desupⁿ)} \isntruncmap{(n-1)}␣f\]\end{theorem}

\begin{proof}

First we note that by Lemma \ref{lma:equiv-Tn-Tinf-homomorphism}, the
type of \(\Tⁿ\)-coalgebra homomorphisms from \((X,m)\) to
\((\Vⁿ_∞,\desupⁿ)\) for which the underlying map is \((n-1)\)-truncated,
is a subtype of the type of \(\T^∞\)-coalgebra homomorphisms from
\((X,m)\) to \((\V^∞_∞,\desup^∞)\). Specifically, we have the following
chain of equivalences and embeddings: \begin{align*}
    ∑_{f : X ↪_{n-1} \Vⁿ_∞} \desupⁿ ∘ f ∼ \Tⁿ␣f ∘ m
        &≃ ∑_{f : X ↪_{n-1} \Vⁿ_∞} \desup^∞ ∘␣π₀ ∘ f ∼ \T^∞␣(π₀ ∘ f) ∘ m \\
        &↪ ∑_{f : X → \Vⁿ_∞} \desup^∞ ∘␣π₀ ∘ f ∼ \T^∞␣(π₀ ∘ f) ∘ m \\
        \label{eq:Tn-hom-subtype-Tinf-hom-step-1}
        &↪ ∑_{f : X → \V^∞_∞} \desup^∞ ∘ f ∼ \T^∞␣f ∘ m
\end{align*} The last step is an instance of the fact that embeddings
are monomorphisms and that dependent sums preserve embeddings.

By Proposition \ref{prop:corec-n-is-Tn-homomorphism}, the first type in
the chain above is inhabited. Since any inhabited type which embeds into
a proposition is contractible, it follows that the first type is
contractible.\end{proof}

Note that this does not contradict Theorem \ref{thm:non-terminality}
since the map \(g\) in the proof of the theorem is not an embedding.

\subsection{The coiterative hierarchy as a model of set
theory}\label{the-coiterative-hierarchy-as-a-model-of-set-theory}

As we recalled in Subsection \ref{fixed-point-models}, Rieger observed
that any fixed point of the powerset functor is a model of ZFC\(^-\)
(ZFC without foundation/regularity) \citep[Theorem III]{rieger1957}. A
corresponding result was shown in Gylterud and Stenholm~\cite{GylterudStenholm2026} for
models of set theory in HoTT---the powerset functor in this case being
\(\T⁰\). Specifically, a fixed point of \(\T⁰\) in HoTT is a model of

\begin{itemize}
\tightlist
\item
  empty set,
\item
  unordered pairing,
\item
  restricted separation,
\item
  replacement,
\item
  union,
\item
  exponentiation,
\item
  infinity/natural numbers.
\end{itemize}

In fact, natural higher type level generalisations of these axioms were
defined, and it was shown that fixed points of \(\Tⁿ\) satisfy the
axioms at level \(n\) or less\footnote{There is also a requirement about
  the fixed point being appropriately locally small.}
\citep[Section~5]{GylterudStenholm2026}. Moreover, the type \(\Vⁿ\) was
shown to be the initial algebra of the functor \(\Tⁿ\) and as such was
shown to model the axiom of foundation, in addition to the ones above.

We shall now look at how \(\Vⁿ_∞\) forms an ∈-structure, and observe
that it satisfies many of the same axioms as \(\Vⁿ\) with the critical
exception of foundation. Instead of foundation, \(\Vⁿ_∞\) will satisfy
Scott's anti-foundation Axiom. The definition of the elementhood
relation on \(\Vⁿ_∞\) is the one which is induced by its coalgebra
structure. The idea is that the elements of a tree are the children of
the root.

\begin{definition}[\agdalink{https://elisabeth.stenholm.one/non-wellfounded-set-theory/v4/e-structure.from-P-inf-coalgebra.html\#1039}]

For \(x,y : \Vⁿ_∞\), define the elementhood relation between them as
\[x ∈ₙ y ≔ \fib␣\widetilde{y}\;x.\]\end{definition}

The relation \(∈ₙ\) is extensional: the canonical map
\[x = y → ∏_{z : \Vⁿ_∞} z ∈ₙ x ≃ z ∈ₙ y\] is an equivalence (Lemma
\ref{coalgebra-structures}). Thus, we have a ∈-structure
\((\Vⁿ_∞,{∈ₙ})\).

The following result is an instance of the results in Section 5 of
Gylterud and Stenholm~\cite{GylterudStenholm2026} (also recalled in Subsection
\ref{fixed-point-models}). The theorem shows that a (locally small)
fixed point of \(\Tⁿ\) models all the defined properties except
foundation. Having shown that \(\Vⁿ_∞\) is such a locally small fixed
point, we apply it to obtain:

\begin{theorem}[\agdalink{https://elisabeth.stenholm.one/non-wellfounded-set-theory/v4/coiterative.set.properties.html}]

\label{thm:Vn-fixed-point-properties}For \(n : \Nat^∞\), \((\Vⁿ_∞,∈ₙ)\)
satisfies the following properties, as defined in
Gylterud and Stenholm~\cite{GylterudStenholm2026}:

\begin{itemize}
\tightlist
\item
  empty set,
\item
  \(U\)-restricted \(n\)-separation,
\item
  ∞-unordered \(I\)-tupling, for all \(k : \Nat_{-1}\) and \(k\)-types
  \(I : U\) such that \(k < n\),
\item
  \(k\)-unordered \(I\)-tupling, for all \(k : \Nat_{-1}\) such that
  \(k ≤ n\) and \(I : U\),
\item
  \(k\)-replacement, for all \(k : \Nat_{-1}\) such that \(k ≤ n\),
\item
  \(k\)-union, for all \(k : \Nat_{-1}\) such that \(k ≤ n\),
\item
  exponentiation, for any ordered pairing structure,
\item
  natural numbers for any \((n-1)\)-truncated representation.
\end{itemize}

\end{theorem}

Since \(\Vⁿ_∞\) is not the initial \(\Tⁿ\)-algebra, it does not follow
that it is a model of foundation. Indeed, since it contains infinite
trees, it is \emph{not} a model of foundation. Neither is \(\Vⁿ_∞\) the
terminal \(\Tⁿ\)-coalgebra, and thus it does not follow that it is a
model of Aczel's anti-foundation axiom. In this section we will show
that it is, however, a model of Scott's anti-foundation axiom.

\subsection{\texorpdfstring{\(\Vⁿ_∞\) models Scott's anti-foundation
axiom}{\textbackslash Vⁿ\_∞ models Scott's anti-foundation axiom}}\label{vux207f_-models-scotts-anti-foundation-axiom}

As \(\Vⁿ_∞\) is not the initial \(\Tⁿ\)-algebra, \((\Vⁿ_∞, ∈ₙ)\) is not
a model of foundation. Indeed, \(\Vⁿ_∞\) contains anti-wellfounded sets,
the simplest one being the infinite unary tree:

\begin{center}
    \begin{tikzcd}
        \bullet \arrow[d] \\
        \bullet \arrow[d] \\
        \vdots
    \end{tikzcd}
\end{center}

As discussed at the start of this paper, there are several
anti-foundation axioms in material set theory. In this section we will
show that \((\Vⁿ_∞, ∈ₙ)\) has Scott \(n\)-anti-foundation.

By Theorem 1 in Gylterud and Stenholm~\cite{GylterudStenholm2026} and Theorem
\ref{thm:Vn-fixed-point-properties}, \((\Vⁿ_∞, ∈ₙ)\) has an ordered
pairing structure. Let \(〈-,-〉 : \Vⁿ_∞ × \Vⁿ_∞ ↪ Vⁿ_∞\) denote this
structure.

\begin{theorem}[\agdalink{https://elisabeth.stenholm.one/non-wellfounded-set-theory/v4/coiterative.set.properties.html\#23590}]

\label{thm:vninf-has-nsafa}For each \(n:ℕ^∞₀\) the ∈-structure
\((\Vⁿ_∞,∈ₙ)\) has the Scott \(k\)-anti-foundation property (\(k\)-SAFA)
for any \(k ≤ n\).\end{theorem}

\begin{proof}

SAFA₂ is immediate from \(\Vⁿ_∞\) being a simple \(\T^∞\)-coalgebra by
Proposition \ref{prop:vnif-simple} and Proposition
\ref{prop:equiv-decoration-Tn-coalgebra}.

For \(k\)-SAFA₁ it suffices to look at the top case \(k = n\). Let
\(g : \Vⁿ_∞\) be a Scott \(n\)-extensional graph. Elements in \(Vⁿ_∞\)
are all \(n\)-types, so by Lemma \ref{lem:truncation-of-n} we have a
\(\Tⁿ\)-coalgebra \(\tar_{n,g} : \Target g → \Tⁿ(\Target g)\). Since
\(g\) is Scott \(n\)-extensional,
\(\corec^∞_{\tar_{n,g}} : \Target g → \V^∞_∞\) is \((n-1)\) truncated
and we obtain by Proposition \ref{prop:corec-lands-in-Vn} and Definition
\ref{corec-n}, a map \(\corec^n : \Target g → Vⁿ_∞\)

To obtain from this a \(\T^∞\)-coalgebra homomorphism from
\((\Vⁿ_∞,m_g)\) to \((\Vⁿ_∞,m_{∈ₙ})\), and thus an ∞-decoration by
Proposition \ref{prop:equiv-decoration-Tn-coalgebra}, let \begin{align*}
   d␣x = \ssupⁿ␣\big(∑_{y:\Vⁿ_∞}〈x,y〉∈g,λ(y,e).\corec^n␣(y,|(x,e)|)\big).
\end{align*} This is a valid application of \(\supⁿ\) since
\(∑_{y:\Vⁿ_∞}〈x,y〉∈g\) is essentially small and \(\corec^n\) is
\((n-1)\)-truncated and thus its composition with the map
\(\left(∑_{y:\Vⁿ_∞}〈x,y〉∈g\right) → \Target␣g\) sending \((y,e)\) to
\((y,|(x,e)|)\) is \((n-1)\)-truncated. It remains to check that the
coalgebra homomorphism square commutes,
i.e.~\(\desup^n␣(d␣x) = \T^∞d␣(m_g␣x)\). Note that the first components
of both \(\desup^n␣(d␣x)\) and \(\T^∞d␣(m_g␣x)\) are
\(∑_{y:\Vⁿ_∞}〈x,y〉∈g\). For the second component we have the following
chain of equalities:
\begin{align}
    π₁␣\left(\T^∞d␣(m_g␣x)\right)
        &= d∘π₀ \label{afa-v-n-1} \\ 
        &= λ(y,e). d␣y  \label{afa-v-n-2}\\
        &= λ(y,e). \ssupⁿ␣\left(∑_{z:\Vⁿ_∞}〈y,z〉∈g,λ(z,e').\corec^n␣(z,|(y,e')|)\right) \label{afa-v-n-3}\\ 
        &= λ(y,e). \ssupⁿ␣\left(\T^∞\corec^n␣\left(∑_{z:\Vⁿ_∞}〈y,z〉∈g,λ(z,e').(z,|(y,e')|)\right)\right) \label{afa-v-n-4}\\ 
        &= λ(y,e). \ssupⁿ␣(\T^∞\corec^n␣(\tar_g␣(y,|(x,e)|))) \label{afa-v-n-5} \\ 
        &= λ(y,e). \corec^n(y,|(x,e)|) \label{afa-v-n-6} \\ 
        &= π₁␣\left(\desup^n␣(d␣x)\right)\label{afa-v-n-7}
\end{align}

Step (\ref{afa-v-n-1}) is the action of \(\T^∞\) and (\ref{afa-v-n-2})
spells out the composition. Writing out the definition of~\(d␣y\) yields
(\ref{afa-v-n-3}). Step (\ref{afa-v-n-4}) is the action of \(\T^∞\) on
morphisms. The definition of \(\tar_g\) yields (\ref{afa-v-n-5}). Step
(\ref{afa-v-n-6}) follows from the fact that \(\corec^n\) is a
\(\T^∞\)-coalgebra homomorphism. And finally, (\ref{afa-v-n-7}) is
obtained by the definition of \(d␣x\).\end{proof}

\section{\texorpdfstring{The terminal
\(\T⁰\)-coalgebra}{The terminal \textbackslash T⁰-coalgebra}}\label{the-terminal-tux2070-coalgebra}

In this section we describe a general construction of terminal
coalgebras for functors satisfying a certain accessibility condition. We
apply this to \(\T⁰\) to obtain a model of Aczel's anti-foundation axiom
in Homotopy Type Theory, assuming propositional resizing. This is a
formalisation in type theory of a theorem due to Aczel and Mendler~\cite{aczel-mendler},
which states that every \emph{set-based} endofunctor on the category of
proper classes has a terminal coalgebra. We describe how to translate
the original proof of Aczel and Mendler, written in the language of set
theory, in Homotopy Type Theory.

Notice that Aczel and Mendler~\cite{aczel-mendler} explicitly assume the axiom of choice
in the paragraph before their Lemma 4.1 and employ the law of excluded
middle in many places, e.g.~in the proof of their Lemma 4.1. Therefore
the translation of their results in the constructive setting of HoTT
requires some care. In the type theoretic statement of the theorem,
proper classes are replaced by large types, and sets are replaced by
small types. The notion of set-based functor is replaced by a certain
accessibility condition with respect to small types. We were able to
remove all invocations of choice principles from the original proof, but
not all impredicativity. In fact, the existence of terminal coalgebras
is guaranteed only under the assumption of
\emph{propositional resizing}, a form of impredicativity for
propositions. Here we recall the principle in a formulation given by
de Jong and Escardó~\cite{JongEscardo23}.

\begin{definition}[\agdalink{https://elisabeth.stenholm.one/non-wellfounded-set-theory/v4/terminal/Utilities.html\#2551}]

The principle of \textbf{propositional resizing} states that every
proposition \(P : \Type\) is essentially small, i.e.~it is equivalent to
a small proposition \(Q : U\).\end{definition}

We do not assume propositional resizing globally, but we mark all
theorems that require its assumption.

Remember that \(\T⁰\) does not have a functorial action on \emph{all}
functions, only on ones with locally small codomain. In the presence of
propositional resizing, these can also be functions with set-valued
codomain. This means that the Aczel--Mendler theorem does not
immediately apply to \(\T⁰\). Nevertheless, in the last part of this
section we will show how to appropriately adjust the statement and proof
of the theorem in order to construct terminal coalgebras also for
``functors'' such as \(\T⁰\).

The constructions in this section are presented using a single universe
\(U : \Type\), in order to keep the presentation consistent with the
previous section, but they can be generalised to polymorphic universes.
More details on the universe-polymorphic constructions can be found in
the Agda formalisation.

\subsection{\texorpdfstring{\(U\)-based
functors}{U-based functors}}\label{u-based-functors}

Aczel and Mendler's theorem applies to set-based endofunctors on proper
classes, where, intuitively, a functor is set-based when its value on a
proper class \(X\) is the colimit of values on small subsets of \(X\).
Before reformulating this accessibility condition in our type theoretic
setting, we recall some definitions and establish some notation.

\textbf{Note:} In this section, we globally assume functors to be
set-valued, i.e \(\F␣X\) is a set, independently of the type level of
\(X\).

\begin{definition}

\label{def:terminal}Let \(A : \Type\) and \(α : A → \F␣A\) be a
coalgebra. We say that \(α\) is

\begin{itemize}
\tightlist
\item
  \textnormal{(}\agdalink{https://elisabeth.stenholm.one/non-wellfounded-set-theory/v4/terminal/Coalgebras.html\#1156}\textnormal{)}
  \textbf{\(U\)-simple} if, for all \(B : U\) and coalgebras
  \(β : B → \F␣B\), the type of coalgebra homomorphisms from \(β\) to
  \(α\) is a proposition;
\item
  \textnormal{(}\agdalink{https://elisabeth.stenholm.one/non-wellfounded-set-theory/v4/terminal/Coalgebras.html\#1449}\textnormal{)}
  \textbf{\(U\)-terminal} if, for all \(B : U\) and coalgebras
  \(β : B → \F␣B\), the type of coalgebra homomorphisms from \(β\) to
  \(α\) is contractible.
\end{itemize}

\end{definition}

Aczel and Mendler write ``strongly extensional'' instead of
``\(U\)-simple''. Assuming propositional resizing, the Aczel--Mendler
theorem guarantees the existence of a \(U\)-terminal coalgebra for every
functor \(\F\). But the existence of a terminal coalgebra is guaranteed
only in case \(\F\) satisfies an accessibility condition. This condition
is a type-theoretic reformulation (and slight generalisation) of Aczel
and Mendler's notion of set-based functor.

\begin{definition}[\agdalink{https://elisabeth.stenholm.one/non-wellfounded-set-theory/v4/terminal/UBased.html\#515}]

\label{def:based}A functor \(\F\) is \textbf{\(U\)-based} if, for any
large type \(X : \Type\) and \(x : \F␣X\), there is a small type
\(Y : U\), a function \(ι : Y → X\) and an element \(y : \F␣Y\) such
that \(\F␣ι\;y = x\).\end{definition}

The existential quantification in the above statement is strong, i.e.~it
is a Σ-type without propositional truncation around it. In other words,
there is a function assigning to each pair \((X : \Type ,x : \F␣X)\) a
tuple \((Y : U,ι : Y → X ,y : \F␣Y,eq: \F␣ι\;y = x)\). Intuitively,
\(\F\) is \(U\)-based when \(\F␣X\) is the colimit of \(\F␣Y\), where
\(Y\) ranges over small generalised elements of \(X\). Notice that the
definition is slightly different from the one of Aczel and Mendler, as
the they require \(Y\) to be a subset of \(X\), i.e.~\(ι\) is an
embedding in their definition. This restriction is not crucial in the
construction of the terminal coalgebra, so we remove it from the
definition.

Notice that Definition \ref{def:based} admits a slight reformulation,
that will become useful later on: a functor is \(U\)-based whenever for
all \(X : \Type\) the function
\[(λ (Y,ι,y).␣\F␣ι\;y) : \left( ∑_{(Y,ι) : \T^∞ X} \F␣Y \right) → \F␣X\]
has a section
\(\operatorname{base}_{\F} : \F␣X → ∑_{(Y,ι) : \T^∞ X} \F␣Y\).

Examples of \(U\)-based functors include all polynomial functors
\(\F␣X ≔ ∑_{a:A} B␣a → X\) with \(A : U\) and \(B : A → U\), i.e.~when
\(A\) and \(B\) are valued in small types. Given \((a,f) : \F␣X\), the
\(U\)-basedness of \(\F\) is evidenced by taking \(Y ≔ B␣a\), \(ι ≔ f\)
and \(y ≔ (a , \id)\). Another canonical example is the functor
\(\T^∞\), whose \(U\)-basedness is evidenced in a similar way, but now
taking \(Y ≔ A\) when given \((A,f) : \T^∞ X\).

\subsection{Relation lifting and
precongruences}\label{relation-lifting-and-precongruences}

There are many ways to lift a (possibly proof-relevant) relation on a
type \(X\) to a relation on \(\F␣X\) \citep{staton-bisim}. Many of these
liftings are well-behaved only when the functor \(\F\) preserves weak
pullbacks. This restriction can be avoided by employing Aczel and
Mendler's notion of relation lifting.

\begin{definition}[\agdalink{https://elisabeth.stenholm.one/non-wellfounded-set-theory/v4/terminal/Utilities.html\#3245}]

Given \(X: \Type\), the \textbf{relation lifting} \(E_{\F}\) takes a
relation \(R : X × X → \Type\) and produces a relation
\(E_{\F}␣R : \F␣X × \F␣X → \Type\) as follows:
\[E_{\F}␣R␣(x,y) ≔ \left(\F␣[-]_R␣x = \F␣[-]_R␣y\right)\] where
\([-]_R\) is the point constructor of the set quotient
\(X/R\).\end{definition}

In HoTT, the set quotient \(X/R\) is defined as a higher inductive type
and the relation \(R\) is not required to be an equivalence relation.
Notice that \(E_{\F}␣R\) is always propositionally-valued since
\(\F␣(X/R)\) is always a set. Even if \(R\) is valued in \(U\) instead
of \(\Type\), there is no guarantee that \(E_{\F}␣R\) is also valued in
\(U\), as \(\F␣(X/R)\) may not be locally \(U\)-small. But this is true
under the assumption of propositional resizing.

\newcommand{\isPrecong}{\ensuremath{\operatorname{is-precong}}}
\newcommand{\Precong}{\ensuremath{\operatorname{Precong}}}

\begin{definition}[\agdalink{https://elisabeth.stenholm.one/non-wellfounded-set-theory/v4/terminal/Precongruences.html\#1142}]

Given a coalgebra \(α : X → \F␣X\), a relation \(R : X × X → \Type\) is
called a \textbf{precongruence} if the following type is inhabited:
\[\ensuremath{\operatorname{is-precong}}_α R ≔ ∏_{x, y : X} R␣(x,y) → E_{\F}␣R␣(α␣x,α␣y)\]\end{definition}

The type of propositionally-valued precongruences on the coalgebra \(α\)
is denoted \(\ensuremath{\operatorname{Precong}}_α\), and we write
\(\ensuremath{\operatorname{Precong}}^U_α\) for the type of
propositionally-valued small precongruences.

\begin{definition}[\agdalink{https://elisabeth.stenholm.one/non-wellfounded-set-theory/v4/terminal/Precongruences.html\#2869}]

A coalgebra \(α : X → \F␣X\) is called \textbf{\(U\)-precongruence
simple} if, for all \(x,y : X\) such that \(R␣(x,y)\) for some reflexive
\(R : \ensuremath{\operatorname{Precong}}^U_α\), then also
\(x = y\).\end{definition}

Aczel and Mendler require the precongruence in the definition of
\(U\)-precongruence simple coalgebra (which they call ``s-extensional'')
to be a congruence, i.e.~an equivalence relation on \(X\). We do not
require symmetry and transitivity, as reflexivity is sufficient for our
purposes (crucially in the proof of Proposition \ref{prop:prec-simple}).
The terminology ``simple'' comes from Rutten~\cite{Rutten2000}, denoting
coalgebras for which bisimulation implies equality. We generalise the
notion from bisimulation to reflexive precongruence.

The greatest (or ``maximal'' in the terminology of
\citep{aczel-mendler}) precongruence on a coalgebra \(α\) is the
propositional truncation of the disjoint union of all its small
precongruences: \begin{equation}\label{eq:greatest-precong}
x ∼_α y ≔ \left\| ∑_{R : \ensuremath{\operatorname{Precong}}^U_α} R␣(x,y) \right\|_{-1}
\end{equation} It is possible to show that
\((∼_α) : \ensuremath{\operatorname{Precong}}_α\).

We can form the set quotient \(X/{∼_α}\), which satisfies a number of
important properties. First, under the assumption of propositional
resizing, \(X/{∼_α}\) is locally \(U\)-small. Second, \(X/{∼_α}\) has an
\(\F\)-coalgebra structure \(α^q : X/{∼_α} → \F␣(X/{∼_α})\) defined by
structural recursion. The case of the point constructor is given as
follows: \(α^q␣[x]_{∼_α} ≔ \F␣[-]_{∼_α} (α␣x)\). The constructor
\([-]_{∼_α}\) is a coalgebra homomorphism between \(α\) and \(α^q\).

\begin{proposition}[\agdalink{https://elisabeth.stenholm.one/non-wellfounded-set-theory/v4/terminal/MaxQuotExt.html\#1472}]

\label{prop:prec-simple}The coalgebra \(α^q : X/{∼_α} → \F␣(X/{∼_α})\)
is \(U\)-precongruence simple.\end{proposition}

\begin{proof}

Applying the elimination principle of set quotients, it is sufficient to
show that given \(x,y:X\), a propositionally-valued reflexive
precongruence \(R : X/{∼_α} × X/{∼_α} → U\) and a proof of
\(R␣([x]_{∼_α},[y]_{∼_α})\), then \(x ∼_α y\). In other words, we need
to find a propositionally-valued precongruence \(S : X × X → U\) such
that \(S␣(x,y)\). Take \(S␣(a,b) ≔ R␣([a]_{∼_α},[b]_{∼_α})\). Notice
that, since the relation \(R\) is reflexive, the types \((X/{∼_α})/R\)
and \(X/S\) are isomorphic, and the underlying function
\(c : (X/{∼_α})/R → X/S\) makes the following square commute:
\begin{equation}\label{eq:iter-quot}
\begin{tikzcd}
    X \arrow[swap, d, "{[-]_{∼_α}}"] \arrow[r, "{[-]_S}"] & {X/S}
\\
{X/{∼_α}} \arrow[swap, r, "{[-]_R}"] & {(X/{∼_α})/R} \arrow[swap, u, "c"]
\end{tikzcd}
\end{equation} Let \(a,b : X\) and suppose \(S␣(a,b)\). The following
sequence of equalities proves that \(S\) is a precongruence:
\begin{align}
\F␣[-]_S␣(α␣a)
\label{terminal-step-1}
&= \F␣(c ∘ [-]_R ∘ [-]_{∼_α})␣(α␣a) \\
\nonumber
&= \F␣c\;(\F␣[-]_R␣(α^q␣[a]_{∼_α})) \\
\label{terminal-step-2}
&= \F␣c\;(\F␣[-]_R␣(α^q␣[b]_{∼_α})) \\
\nonumber
&= \F␣(c ∘ [-]_R ∘ [-]_{∼_α})␣(α␣b) \\
\label{terminal-step-3}
&= \F␣[-]_S␣(α␣b)
\end{align} Step (\ref{terminal-step-1}) follows by (\ref{eq:iter-quot})
and step (\ref{terminal-step-2}) is the fact that R is a precongruence.
Finally, in step (\ref{terminal-step-3}) we use (\ref{eq:iter-quot})
again.\end{proof}

\begin{proposition}[\agdalink{https://elisabeth.stenholm.one/non-wellfounded-set-theory/v4/terminal/Precongruences.html\#4812}]

\label{prop:precongsimp-to-simp}Every \(U\)-precongruence simple
coalgebra with locally \(U\)-small carrier is
\(U\)-simple.\end{proposition}

\begin{proof}

Let \(α : X → \F␣X\) be a \(U\)-precongruence simple coalgebra with
\(X\) locally \(U\)-small. Let \(f,g\) be two coalgebra homomorphisms
from another coalgebra \(β : Y → \F␣Y\) to \(α\). Given \(y : Y\), it is
sufficient to show that \(f␣y = g␣y\) (remember that we globally assume
the functor \(\F\) to be set-valued). From the precongruence simplicity
of \(α\), it is sufficient to find a propositionally-valued reflexive
precongruence \(S : X × X → U\) relating \(f␣y\) and \(g␣y\). Consider
the relation: \[R'␣x␣x' ≔ ∑_{y : Y} (x = f␣y) × (x' = g␣y)\] and its
propositional reflexive closure
\(R␣x␣x' ≔ \left\| R'␣x␣x' + (x = x')\right\|_{-1}\). It is not hard to
show that \(R\) is a precongruence on \(α\), which moreover relates
\(f␣y\) and \(g␣y\) as \(| \inl (y , \refl , \refl) | : R␣(f␣y)␣(g␣y)\).

If \(X\) is a large type then the relation \(R\) is valued in \(\Type\).
But since \(Y : U\) and \(X\) is locally \(U\)-small, there is a
\(U\)-valued relation \(S : X × X \to U\) such that
\(S␣x␣x' \simeq R␣x␣x'\) for all \(x,x' : X\). Moreover,
\(S : \ensuremath{\operatorname{Precong}}^U_α\) and \(S\) relates
\(f␣y\) and \(g␣y\).\end{proof}

\begin{corollary}

\label{cor:str-ext}Assuming propositional resizing, the coalgebra
\(α^q : X/{∼_α} → \F␣(X/{∼_α})\) is \(U\)-simple.\end{corollary}

\begin{proof}

The coalgebra \(α^q\) is \(U\)-precongruence simple by Proposition
\ref{prop:prec-simple}. \(X/{∼_α}\) is locally \(U\)-small by
propositional resizing. Therefore \(α^q\) is \(U\)-simple by Proposition
\ref{prop:precongsimp-to-simp}.\end{proof}

\subsection{\texorpdfstring{The \(U\)-terminal
coalgebra}{The U-terminal coalgebra}}\label{the-u-terminal-coalgebra}

The \(U\)-terminal coalgebra of a functor \(\F\) is built in two steps.
First, define the \emph{weakly} \(U\)-terminal coalgebra as the disjoint
union of all small coalgebras: \begin{equation}\label{eq:weakly}
wν\F_U ≔ ∑_{X : U}∑_{α : X → \F␣X} X.
\end{equation} Every small coalgebra \(α : X → \F␣X\) clearly injects in
the union \(α^* : X → wν\F_U\), \(α^*␣x ≔ (X,α,x)\). The coalgebra
structure \(ζ : wν\F_U → \F␣(wν\F_U)\) is given by
\(ζ␣(X,α,x) ≔ \F␣α^*␣(α␣x)\). It easy to prove that \(α^*\) is a
coalgebra homomorphism between \(α\) and \(ζ\).

In order to turn the weakly \(U\)-terminal coalgebra into a
\emph{strong} \(U\)-terminal coalgebra, we quotient its carrier
\(wν\F_U\) by the greatest precongruence on \(ζ\) (introduced in
(\ref{eq:greatest-precong})): \(ν\F_U ≔ wν\F_U/{∼_{ζ}}\). We know this
has a coalgebra structure \(ζ^q\). Moreover, given a small coalgebra
\(α : X → \F␣X\), there is a coalgebra homomorphism from it to \(ζ^q\)
given by the composition of \(α^*\) and \([-]_{∼_ζ}\). Invoking
Corollary \ref{cor:str-ext}, which assumes propositional resizing, we
know that this is the only such coalgebra homomorphism.

\begin{theorem}[\agdalink{https://elisabeth.stenholm.one/non-wellfounded-set-theory/v4/terminal/UTerminal.html\#4265}]

\label{thm:aczel-mendler1}Assuming propositional resizing, the coalgebra
\(ζ^q : ν\F_U → \F␣(ν\F_U)\) is \(U\)-terminal.\end{theorem}

\subsection{\texorpdfstring{\label{sec:aczel-mendler}The Aczel--Mendler
theorem}{The Aczel--Mendler theorem}}\label{the-aczelmendler-theorem}

We finally show how the \(U\)-terminal coalgebra \(ζ^q\) is also
\emph{terminal} with respect to large coalgebras, provided the functor
\(\F\) is \(U\)-based.

\newcommand{\bind}{\ensuremath{\operatorname{bind}}}

First, notice that \(\T^∞\) is not only a polynomial functor, but a
polynomial monad. Its unit \(η : X → \T^∞ X\) is \(η␣x ≔ (1,λ *.x)\).
The Kleisli extension
\(\ensuremath{\operatorname{bind}}g: \T^∞ X → \T^∞ Y\) of a function
\(g : X → \T^∞ Y\) is obtained by forming the disjoint union of all
indexing types:
\[\ensuremath{\operatorname{bind}}g␣(A,f) ≔ \left(∑_{a : A} π_0␣(g (f a)),␣λ(a,y).␣π_1␣(g (f a)) y
             \right)\] Given \(g : X → \T^∞ X\), its Kleisli extension
can be iterated a finite number of times:
\begin{equation*}\label{eq:iter-bind}
\begin{array}{l@{}l}
\multicolumn{2}{l}{\ensuremath{\operatorname{bind}}: \mathbb{N} → (X → \T^∞ X) → \T^∞ X → \T^∞ X}\\ 
\ensuremath{\operatorname{bind}}^0 &g␣z ≔ z \\
\ensuremath{\operatorname{bind}}^{n+1} &g␣z ≔ \ensuremath{\operatorname{bind}}g␣(\ensuremath{\operatorname{bind}}^n g␣z).
\end{array}
\end{equation*} It can also be iterated an infinite number of times, by
collecting all the finite approximations: \begin{align*}
&\ensuremath{\operatorname{bind}}^∞ g: \T^∞ X → \T^∞ X \\ 
&\ensuremath{\operatorname{bind}}^∞ g␣z ≔ \left(∑_{n : ℕ} π_0␣(\ensuremath{\operatorname{bind}}^n g␣z),␣λ(n,x).␣π_1␣(\ensuremath{\operatorname{bind}}^n g␣z) x \right)
\end{align*}

Given a large coalgebra \(α : X → \F␣X\) for a \(U\)-based functor
\(\F\), one can construct a \(\T^∞\)-coalgebra structure on \(X\) as
follows: \(\widehat{α}␣x ≔ π_0␣(\operatorname{base}_{\F}␣(α␣x))\).

\begin{proposition}[\agdalink{https://elisabeth.stenholm.one/non-wellfounded-set-theory/v4/terminal/UBased.html\#1895}]

\label{prop:onestep}Let \(\F\) be a \(U\)-based functor and
\(α : X → \F␣X\) a large coalgebra. For all \(z : \T^∞ X\), there is a
function
\(α_z : π₀␣z → \F␣(π₀␣(\ensuremath{\operatorname{bind}}\widehat{α}␣z))\)
such that the following diagram commutes:

\begin{center}
\begin{tikzcd}
    {π₀␣z} \arrow[swap, d, "α_z"] \arrow[rrr, "π₁␣z"] & & & X \arrow[d, "α"] \\
    {\F␣(π₀␣(\bind \widehat{α}␣z))} \arrow[swap, rrr, "\F␣(π₁␣(\bind \widehat{α}␣z))"] & & & \F␣X                 
\end{tikzcd}
\end{center}

\end{proposition}

\begin{proof}

Let \(a : π₀␣z\). Since \(\F\) is \(U\)-based, there exist \(A : U\),
\(ι : A → X\) and \(y : \F␣A\) such that \(\F␣ι\;y = α␣(π₁␣z␣a)\). In
other words \(y ≡ π₂␣(\operatorname{base}_{\F}␣(α␣(π₁␣z␣a)))\). Take
\(α_z␣a ≔ \F␣(λx.␣(a,x))␣y\).\end{proof}

The construction of Proposition \ref{prop:onestep} can be iterated,
producing a family of functions
\[α^n_z : π₀␣(\ensuremath{\operatorname{bind}}^n \widehat{α}␣z) → \F␣(π₀␣(\ensuremath{\operatorname{bind}}^{n+1} \widehat{α}␣z))\]
indexed by a natural number \(n\), which makes the following family of
diagrams commute: \begin{equation}\label{eq:diagram}
    \begin{tikzcd}
        {π₀␣(\ensuremath{\operatorname{bind}}^n \widehat{α}␣z)} \arrow[swap, d, "α^n_z"] \arrow[rrrr, "π₁␣(\ensuremath{\operatorname{bind}}^n \widehat{α}␣z)"] & & & & X \arrow[d, "α"] \\
        {\F␣(π₀␣(\ensuremath{\operatorname{bind}}^{n+1} \widehat{α}␣z))} \arrow[swap, rrrr, "\F␣(π₁␣(\ensuremath{\operatorname{bind}}^{n+1} \widehat{α}␣z))"] & & & & \F␣X                 
    \end{tikzcd}
\end{equation}

\begin{proposition}[\agdalink{https://elisabeth.stenholm.one/non-wellfounded-set-theory/v4/terminal/UBased.html\#2585}]

\label{prop:small-coalgebra}Let \(\F\) be a \(U\)-based functor and
\(α : X → \F␣X\) a large coalgebra. Then each \(z : \T^∞ X\) determines
a small coalgebra \(α^∞_z : X_z → \F␣(X_z)\) and a coalgebra
homomorphism \(k_z\) from \(α^∞_z\) to \(α\).\end{proposition}

\begin{proof}

Define the carrier \(X_z\) as
\(π₀␣(\ensuremath{\operatorname{bind}}^∞ \widehat{α}␣z)\) and its
coalgebra structure as \[α^∞_z␣(n,x) ≔ \F␣(λy.␣n+1,y)␣(α^n_z␣x).\] There
is a function
\(k_z␣(n,x) ≔ π₁␣(\ensuremath{\operatorname{bind}}^n \widehat{α}␣z)␣x\)
between \(X_z\) and \(X\). The fact that this is a coalgebra
homomorphism between \(α^∞_z\) and \(α\) follows from the commutativity
of the family of diagrams in (\ref{eq:diagram}).\end{proof}

Notice also the existence of a function \(u_z : π₀␣z → X_z\) sending
\(x\) to the pair \((0,x)\), which makes the triangle below commute.
Since \(k_z\) is a coalgebra homomorphism, the square below also
commutes: \begin{equation}\label{eq:up}
    \begin{tikzcd}
        & {π₀␣z} \arrow[swap, dl, "u_z"] \arrow[dr, "π₁␣z"] & \\
    {X_z} \arrow[swap, d, "α^∞_z"] \arrow[rr, "k_z"] & & X \arrow[d, "α"] \\
        {\F␣X_z} \arrow[swap, rr, "\F␣k_z"] & & \F␣X                 
    \end{tikzcd}
\end{equation}

Given \(z : \T^∞ X\) and \(w : \T^∞ X_z\), the latter also determines an
element \(w' : \T^∞ X\) as follows: \(w' ≔ \T^∞ k_z␣w\). The small
coalgebras associated to \(z\) and \(w'\) by Proposition
\ref{prop:small-coalgebra} are in a strong relationship with each other.

\begin{lemma}[\agdalink{https://elisabeth.stenholm.one/non-wellfounded-set-theory/v4/terminal/UBased.html\#3436}]

\label{lemma:submultiset}Let \(\F\) be a \(U\)-based functor and
\(α : X → \F␣X\) a large coalgebra. For all \(z : \T^∞ X\) and
\(w : \T^∞ X_z\), there is a coalgebra homomorphism \(l_{z,w}\) between
\(α^∞_{w'}\) and \(α^∞_z\) that makes the following diagram commute:
\begin{equation}\label{eq:lemma}
\begin{tikzcd}
π₀␣w' \arrow[swap, d, "u_{w'}"] \arrow[equal, rr] & & π₀ w\arrow[d, "π₁␣w"] \\
    {X_{w'}}  \arrow[swap, rr, "l_{z,w}"] & & {X_z} 
\end{tikzcd}
\end{equation}\end{lemma}

\begin{proof}

We only sketch the construction of \(l_{z,w}\). Its definition follows
from the construction of a term \begin{equation*}
  l'_{z,w} : ∏_{(n,x) : X_{w'}} ∑_{(m,y) : X_z} π₁␣(\ensuremath{\operatorname{bind}}^n \widehat{α}␣w')␣x = π₁␣(\ensuremath{\operatorname{bind}}^m \widehat{α}␣z)␣y
  \end{equation*} by taking \(l_{z,w} ≔ π₀ ∘ l'_{z,w}\). The term
\(l'_{z,w}␣(n,x)\) is defined by induction on \(n\). If \(n = 0\), we
return \((π₁␣w␣x,\refl)\). If \(n = n' + 1\), then \(x\) is a pair
\((x',f)\) consisting of
\(x' : π₀␣(\ensuremath{\operatorname{bind}}^{n'}␣\widehat{α}␣w')\) and
\(f : π₀␣(\widehat{α}␣(π₁␣(\ensuremath{\operatorname{bind}}^{n'}␣\widehat{α}␣w')␣x'))\).
In particular, \((n',x') : X_{w'}\). The recursive call
\(l'_{z,w}␣(n',x')\) gives us a tuple consisting of a natural number
\(m : ℕ\), a term
\(y : π₀␣(\ensuremath{\operatorname{bind}}^m␣\widehat{α}␣z)\) and an
equality proof
\(eq : π₁␣(\ensuremath{\operatorname{bind}}^{n'} \widehat{α}␣w')␣x' = π₁␣(\ensuremath{\operatorname{bind}}^m \widehat{α}␣z)␣y\).
We return the tuple \((m +1,(y,f')) : X_{z}\) where
\(f' : π₀␣(\widehat{α}␣(π₁␣(\ensuremath{\operatorname{bind}}^m␣\widehat{α}␣z)␣y))\)
is obtained from \(f\) by transporting along \(eq\).\end{proof}

We are now ready to prove the main result of Aczel and Mendler~\cite{aczel-mendler}.

\begin{theorem}[\agdalink{https://elisabeth.stenholm.one/non-wellfounded-set-theory/v4/terminal/Terminal.html\#634}]

\label{thm:aczel-mendler2}For a \(U\)-based functor, any \(U\)-terminal
coalgebra is also terminal.\end{theorem}

\begin{proof}

Let \(β : Y → \F␣Y\) be a \(U\)-terminal coalgebra and let
\(α : X → \F␣X\) be a large coalgebra. We construct a coalgebra
homomorphism from \(α\) to \(β\). Given \(x : X\), we get
\(η␣x : \T^∞ X\) and therefore, by Proposition
\ref{prop:small-coalgebra}, a small coalgebra
\(α^∞_{η␣x} : X_{η␣x} → \F␣(X_{η␣x})\). From \(U\)-terminality, there
exists a unique coalgebra homomorphism \(h_x\) between \(α^∞_{η␣x}\) and
\(β\).

We now show how this homomorphism can be lifted to one initiating from
the large coalgebra \(α\). First, a function \(h : X → Y\) can be
defined as \(h␣x ≔ h_x␣(u_{η␣x}␣*)\), which is a coalgebra homomorphism.
In order to show this, we need to prove another equation
\(h ∘ k_{η␣x} = h_x\). Let \(a : X_{η␣x}\) and define \(a' : X\) as
\(a' ≔ k_{η␣x}␣a\). We have the following sequence of equalities:
\begin{equation}\label{eq:extra}
h␣(k_{η␣x}␣a)
≡
h_{a'}␣(u_{η␣a'}␣*)
=
h_x␣(l_{η␣x,η␣a}␣(u_{η␣a'}␣*))
=
h_x␣(π₁␣(η␣a)␣*)
≡
h_x␣a
\end{equation} The second equality holds since \(h_{a'}\) is the unique
coalgebra homomorphism from \(α^∞_{η␣a'}\) to \(β\), and the fact that
\(h_x\) and \(l_{η␣x,η␣a}\) (which was introduced in Lemma
\ref{lemma:submultiset}) are both coalgebra homomorphisms. The third
equality is an instance of (\ref{eq:lemma}).

Proving that \(h : X → Y\) is a coalgebra homomorphism is evidenced by
the following sequence of equations, where step (\ref{terminal-step-4})
follows from (\ref{eq:up}), step (\ref{terminal-step-6}) follows from
(\ref{eq:extra}) and step (\ref{terminal-step-5}) is the fact that
\(h_x\) is a coalgebra homomorphism. \begin{align}
\F␣h\;(α␣x)
\label{terminal-step-4}
&= \F␣h\;(\F␣k_{η␣x}␣(α^∞_{η␣x}␣(u_{η␣x}␣*))) \\
\nonumber
&= \F␣(h ∘ k_{η␣x})␣(α^∞_{η␣x}␣(u_{η␣x}␣*)) \\
\label{terminal-step-6}
&= \F␣h_x␣(α^∞_{η␣x}␣(u_{η␣x}␣*)) \\
\label{terminal-step-5}
&= β␣(h_x␣(u_{η␣x}␣*)) \\
\nonumber
&≡ β␣(h␣x)
\end{align}

The coalgebra homomorphism \(h\) is unique. Given another one \(h'\) and
an element \(x : X\), we have the following sequence of equalities: \[
h␣x
≡
h_x␣(u_{η␣x}␣*)
=
h'␣(k_{η␣x}␣(u_{η␣x}␣*))
=
h'␣(π₁␣(η␣x)␣*)
≡
h'␣x
\] The second equality holds since \(h_{x}\) is the unique coalgebra
homomorphism from \(α^∞_{η␣x}\) to \(β\), and the fact that \(h'\) and
\(k_{η␣x}\) are both coalgebra homomorphisms. The third equality is an
instance of the triangle in (\ref{eq:up}).\end{proof}

Putting together Theorems \ref{thm:aczel-mendler1} and
\ref{thm:aczel-mendler2}, we obtain the general terminal coalgebra
theorem of Aczel and Mendler. Assuming propositional resizing, there is
a \(U\)-terminal coalgebra \(ζ^q : ν\F_U → \F␣(ν\F_U)\) for any functor
\(\F\). If the latter happens to be \(U\)-based, then this coalgebra is
also terminal with respect to large coalgebras.

\begin{theorem}[\agdalink{https://elisabeth.stenholm.one/non-wellfounded-set-theory/v4/terminal/Terminal.html\#3551}]

\label{thm:aczel-mendler}Let \(\F\) be a \(U\)-based functor. Assuming
propositional resizing, the coalgebra \(ζ^q : ν\F_U → \F␣(ν\F_U)\) is
terminal.\end{theorem}

\subsection{\texorpdfstring{Adjusting the theorem for
\(\T⁰\)}{Adjusting the theorem for \textbackslash T⁰}}\label{adjusting-the-theorem-for-tux2070}

The powerset construction \(\T⁰\) is not a functor, as it only acts on
functions \(f : X → Y\) with locally small codomain. The type \(Y\) can
also be restricted to be a set if one assumes propositional resizing.
Crucially this means that the Aczel--Mendler theorem described so far
does not apply to it. Luckily, this can be remedied with a few small
modifications.

First, let us call \(\F\) a \textbf{set-valued functor} if \(\F␣X\) is a
set and \(\F\) acts exclusively on set-valued functions, i.e.~its action
on functions is typed \(∏_{X:\Type, Y : \Set} (X → Y) → \F␣X → \F␣Y\).
Clearly \(\T⁰\) is a set-valued functor in this sense, assuming
propositional resizing.

The notion of \(U\)-basedness in Definition \ref{def:based} also needs
to be adjusted. Let \(\Set_U\) be the type of sets in \(U\). We now say
that a set-valued functor is \textbf{\(\Set_U\)-based} if, for any large
\emph{set} \(X : \Set\) and \(x : \F␣X\), there is a small \emph{set}
\(Y : \Set_U\), a function \(ι : Y → X\) and element \(y : \F
Y\) such that \(\F␣ι\;y = x\). In other words, both \(X\) and \(Y\) in
the definition are required to be sets. This is important for the
results of Section \ref{sec:aczel-mendler} to go through when functors
only act on set-valued functions. For example, the bottom functions in
(\ref{eq:diagram}) and (\ref{eq:up}) are well-defined only if \(X\) is a
set. Similarly, the functions \(l_{z,w}\) in Lemma
\ref{lemma:submultiset} can only be coalgebra morphisms in case \(X_z\)
is a set.

\begin{proposition}[\agdalink{https://elisabeth.stenholm.one/non-wellfounded-set-theory/v4/terminal-var/Powerset.html\#5853}]

\label{prop:P-based}\(\T⁰\) is \(\Set_U\)-based.\end{proposition}

\begin{proof}

Let \(X : \Set\) and \(x : \T⁰ X\). Notice that \(π₀␣x : U\) is a set,
since \(π₁␣x : π₀␣x → X\) is an embedding and \(X\) is a set. Therefore
we can return the triple consisting of the small set \(π₀␣x\), the
function \(π₁␣x : π₀␣x → X\) and the element
\((π₀␣x,\id) : \T⁰ (π₀␣x)\).\end{proof}

The weakly \(U\)-terminal coalgebra in (\ref{eq:weakly}) also needs to
be modified. This is because \(wν\F_U\) is not a set, so there cannot be
any coalgebra homomorphism targeting it. The solution is to take its
\emph{set truncation} \(\left\| wν\F_U \right\|_0\) instead. It is
straightforward to define a coalgebra structure on it using the
elimination principle of set truncation.

Finally, assuming that \(X\) is a set in the definition of
\(\Set_U\)-basedness restricts the notion of terminal coalgebra in
Definition \ref{def:terminal} to work only for coalgebras with a set
carrier. We say that a coalgebra \(α : A → \F␣A\) is \textbf{terminal
with respect to sets} if, for all \(B : \Set\) and coalgebras
\(β : B → \F␣B\), the type of coalgebra homomorphisms from \(β\) to
\(α\) is contractible.

With all these restrictions in place, the Aczel--Mendler Theorem
\ref{thm:aczel-mendler} still works.

\begin{theorem}[\agdalink{https://elisabeth.stenholm.one/non-wellfounded-set-theory/v4/terminal-var/SetTerminal.html\#3731}]

Let \(\F\) be a \(\Set_U\)-based set-valued functor. Assuming
propositional resizing, the coalgebra \(ζ^q : ν\F_U → \F␣(ν\F_U)\) is
terminal with respect to sets.\end{theorem}

Since \(ν\F_U\) is itself a set, the theorem implies that there is a
unique coalgebra morphism from \(ζ^q\) to itself, given by the identity
function.

As a corollary, we obtain a terminal coalgebra with respect to sets for
the powerset functor \(\T⁰\). As such, the latter coalgebra validates
Aczel's anti-foundation axiom, by Theorem
\ref{thm:AFA-from-terminal-coalg} and the remark following it.

\begin{corollary}[\agdalink{https://elisabeth.stenholm.one/non-wellfounded-set-theory/v4/terminal-var/Powerset.html\#6701}]

\label{cor:p0-has-terminal-coalgebra}Assuming propositional resizing,
\(\T⁰\) admits a terminal coalgebra with respect to sets. This terminal
coalgebra forms a model of Aczel's anti-foundation axiom.\end{corollary}

\section{Conclusion and future work}\label{conclusion-and-future-work}

In this paper we constructed a non-initial and non-terminal fixed point
of the (restricted) powerset functor and showed that it is a model of
material set theory with Scott's anti-foundation axiom. Moreover, we
constructed the terminal coalgebra of the same functor, assuming
propositional resizing. This is then a model of material set theory with
Aczel's anti-foundation axiom.

There are still questions that remain unanswered, especially the initial
motivation of this paper: to construct the terminal coalgebra of the
powerset functor. The construction in the last section relies in a
crucial way on propositional resizing. Is there a way to construct the
terminal coalgebra, without this assumption? Is it possible to show that
assuming the existence of the terminal coalgebra implies some classical
principle? Or is it independent altogether? Also, given propositional
resizing, is there an easy way to extend Aczel--Mendler to higher type
levels? Would such a result lead to a terminal algebra for \(\Tⁿ\) for
all \(n\)?

\section*{Acknowledgments}
Niccolò Veltri was supported by the Estonian Research Council grant
PSG749.

\bibliographystyle{alphaurl}
\bibliography{list}

\end{document}